\documentclass[12pt]{article}
\usepackage{amssymb}
\newtheorem{theorem}{Theorem}

\newtheorem{proposition}[theorem]{Proposition}
\newtheorem{lemma}[theorem]{Lemma}
\newtheorem{corollary}[theorem]{Corollary}
\newtheorem{remark}[theorem]{Remark}
\newtheorem{remarks}[theorem]{Remarks}

\newcommand{\R}{\mathbb{R}}
\newcommand{\Q}{\mathbb{Q}}

\newcommand{\Sf}{\mathbb{S}}

\newcommand{\Les}{\mathbb{L}}
\newcommand{\Ve}{\mathbb{V}}

\newcommand{\Hy}{\mathbb{H}}
\newcommand{\Oes}{\mathbb{O}}

\newcommand{\spa}{\mbox{span}}
\newcommand{\hess}{\mbox{Hess\,}}

\newcommand{\Tes}{\mbox{T}}

\newcommand{\Ee}{\mathbb {E}}
\newcommand{\grad}{\mbox{grad}}

\newcommand{\po}{{\hspace*{-1ex}}{\bf .  }}
\newcommand{\ii}{isometric immersion }
\newcommand{\iis}{isometric immersions }
\newcommand{\su}{submanifold }
\newcommand{\sus}{submanifolds }
\newcommand{\fnb}{flat normal bundle }
\newcommand{\sffs}{second fundamental forms }
\newcommand{\sff}{second fundamental form }
\newcommand{\rt}{Ribaucour transformation }

\newcommand{\rtf}{Ribaucour transform }

\newcommand\fall{\;\;\mbox{for all}\;\;}
\def\P{{\cal P}}

\def\Fes{{\cal F}}

\def\<{\langle}

\def\>{\rangle}
\def\a{\alpha}
\def\va{\varphi}

\def\bea{\begin{eqnarray*} }
\def\eea{\end{eqnarray*} }
\def\be{\begin{equation} }
\def\ee{\end{equation} }

\def\nap{\nabla^\perp}
\def\proof{\noindent{\it Proof: }}
\def\qed{\ifhmode\unskip\nobreak\fi\ifmmode\ifinner
\else\hskip5 pt \fi\fi\hbox{\hskip5 pt \vrule width4 pt
height6 pt  depth1.5 pt \hskip 1pt }}
\begin{document}

\title{Isothermic submanifolds of Euclidean space}
\author { Ruy Tojeiro}
\date{}
\maketitle

\noindent {\bf Abstract:} {\small We study the problem posed by F. Burstall of developing a theory of isothermic Euclidean submanifolds of dimension greater than or equal to three. As a natural extension of the definition in the surface case, we  call a Euclidean submanifold {\it isothermic\/} if it is locally the image of a conformal immersion of a Riemannian product of Riemannian manifolds whose second fundamental form is adapted to the product net of the manifold. Our main result is a complete classification of all such conformal immersions of Riemannian products of dimension greater than or equal to three. We derive several consequences of this result. We also study whether the classical characterizations of isothermic surfaces as solutions of Christoffel's problem and as envelopes of  nontrivial conformal sphere congruences extend to higher dimensions. } \vspace{2ex}

\vspace{2ex}


\section[Introduction]{
Introduction}

     The classical theory of isothermic surfaces has received a recent renaissance in interest mostly due to
     its connection with the modern theory of integrable systems
     (cf. \cite{bu} and the references therein).
     Recall that a surface in Euclidean space $\R^3$ is {\em isothermic\/} if, away from
     umbilic points,
     its curvature lines form an isothermic net, that is,  there exist locally conformal
     coordinates that diagonalize  the second fundamental form. Standard examples include
     cylinders, cones, surfaces of revolution, quadrics,  surfaces of constant mean curvature
     and, the concept being a conformally invariant one,  the images of such surfaces by an
     inversion in $\R^3$.

     Isothermic surfaces have a number of characterizations from various points of view.
     One of them comes with the problem posed by Christoffel  \cite{chris}
     of determining the surfaces $f\colon\;M^2\to \R^3$ that admit
     locally a {\em dual\/} surface $\Fes\colon\;M^2\to \R^3$ with parallel tangent planes
     to those of $f$
     which induces the same conformal structure but opposite orientation on $M^2$.
     It was shown by Christoffel himself that this property characterizes isothermic
     surfaces, and dual isothermic surfaces are now said to be {\em Christoffel transforms\/}
     one of each other.

  Another characterization involves the notion of {\em Ribaucour transformation\/}.
  Classically, two surfaces in $\R^3$  are said to be related by a \rt when they are
  the enveloping surfaces of a two-parameter sphere congruence such that their curvature
  lines correspond.
  It was proved by Darboux (\cite{da1}) that if two surfaces in $\R^3$ are related
  by a {\em conformal\/} Ribaucour transformation that reverts orientation then they must be
  a {\em Darboux pair\/} of isothermic surfaces.

A nice account of classical aspects of the subject as well as from
the perspective of the modern theory of integrable systems and
loop groups was given in \cite{bu}. In the concluding remarks of
that paper, a list of open problems is posed that suggests
interesting guides for future research. We quote here the first of
them:\vspace{1ex}

{\em Is there any interesting theory of isothermic submanifolds of
$\R^n$ of dimension greater than two? The problem here is to find
a suitable definition that is not too restrictive...  ...One way
forward might be to study submanifolds admitting a conformal
sphere congruence.}\vspace{1ex}

   This article grew out as an attempt to give an answer to this question.
   The preliminary step of defining higher dimensional analogues of isothermic nets
    and  geometrically characterizing them was considered in \cite{to}. This turned out
    to fit into the general problem of studying local decomposability properties of {\em
  orthogonal nets\/} on Riemannian manifolds.  Recall that an orthogonal $k$-net
    ${\cal E}=(E_i)_{i=1,\ldots,k}$ on a Riemannian manifold  $M$
    is a splitting $TM=\oplus_{i=1}^k E_i$ of the tangent bundle of $M$
    by a family of mutually orthogonal integrable subbundles. It is natural to require a
     higher dimensional isothermic
    net ${\cal E}=(E_i)_{i=1,\ldots,k}$ on a Riemannian manifold $M$ to have the  property
      that each point
    $p\in M$ lies in an open neighborhood $U$ such that there exists
    a conformal diffeomorphism of a Riemannian product $\Pi_{i=1}^k M_i$ onto $U$ that maps
    the leaves of the product foliation induced by $M_i$ onto the
    leaves of  $E_i$, $1\leq i\leq k$. A geometric characterization of
    such orthogonal nets was given in Theorem $4.3$ of \cite{to}
    (see Theorem \ref{cor:drh} below), which may be regarded as a conformal
    version of the local decomposition theorem of de Rham. Following a more
    suggestive terminology used in \cite{mrs} for other types of orthogonal nets,
    they were named  {\em conformal product nets\/}
    in \cite{to}, or $CP$-nets for short.

    The next step is to look for the appropriate extension of the property that
    {\em the curvature lines\/}
    form an isothermic net. A little thought then makes it natural to call  a Euclidean
    submanifold {\em isothermic\/}, or more precisely,
    {\em $k$-isothermic\/},
    if it carries a k-$CP$-net ${\cal E}=(E_i)_{i=1,\ldots,k}$ to which its \sff is adapted,
    in the sense that
    each subbundle $E_i$ is invariant by all shape operators of the submanifold.

    Summing things up,   a submanifold is  isothermic if it is locally given by a conformal
    immersion of a Riemannian product whose  second fundamental form  is adapted
    to the product net of the manifold.
    Thus, in studying isothermic Euclidean submanifolds of dimension $n\geq 3$  we are faced
    with the following  problem:\vspace{.5ex}

    {\em Describe the conformal immersions of a Riemannian product of
dimension $n\geq 3$ into Euclidean space whose second fundamental
forms are adapted to its product net.}\vspace{.5ex}

 The isometric version of this problem (in any dimension $n\geq 2$) has led to a useful result
 due to Moore (\cite{mo}), who proved that an isometric immersion of a Riemannian product into
 Euclidean space whose second fundamental form is adapted to the product net of the manifold
 must be an extrinsic
   product of isometric immersions of the factors. Notice that the composition of such an
   extrinsic product
    of $k$ isometric immersions with a conformal transformation of Euclidean space provides
    trivial examples
    of $k$-isothermic submanifolds for any value of $k$. These examples can be seen as the
    higher
    dimensional analogues of the images by an inversion of cylinders in $\R^3$ over plane curves.
    Similar examples can be produced by composing an extrinsic  product of isometric immersions
    into a sphere  with a stereographic projection of that sphere onto Euclidean space.

 Further examples of isothermic submanifolds arise by means of a conformal diffeomorphism
      $\Theta\colon\,\Hy^m(-c)\times \Sf^{N-m}(c)\to \R^N$, $1\leq m\leq N-1$, onto
      the complement of an $(m\!-\!1)$-dimensional sphere  (a point if $m=1$)
      in  $\R^N$ (see Section $3$).
Here $\Hy^m(-c)=\{X\in \Les^{m+1}\!: \<X,X\>\!=-1/c\}$ and
\mbox{$\Sf^{N-m}(c)=\{X\in \R^{N-m+1}\!: \<X,X\>=1/c\}$,} where
$\<\;,\;\>$ stands for the Lorentz inner product in Minkowski
space $\Les^{m+1}$  in the first case and the standard Euclidean
inner product in the latter. Thus, if  $m\geq 2$ (resp., $N-m\geq
2$) then $\Hy^m(-c)$ (resp., $\Sf^{N-m}(c)$) is the standard model
of hyperbolic space (resp., the sphere) with constant sectional
curvature $-c$ (resp., $c$) as an umbilical hypersurface of
$\Les^{m+1}$ (resp., $\R^{N-m+1}$).  Then, for any  \ii
$f\colon\,M_1\to \Hy^m(-c)$ and any extrinsic product
      $g\colon\,M_2\times \ldots\times M_k\to \Sf^{N-m}(c)$ of isometric immersions,
       the conformal immersion $\Theta\circ (f\times g)$ defines a new example of a
       $k$-isothermic submanifold.
      In two dimensions, these examples reduce to conformal
      parameterizations by curvature lines
      of the images of cones and surfaces of revolution  by an inversion.

      The abundance of isothermic surfaces in $\R^3$ gives a hint that their higher dimensional
      analogues should be equally as plentiful.  The hint turns out to be  surprisingly
      misleading: we are able to show that the examples just described  actually comprise
      all isothermic submanifolds of dimension $n\geq 3$.
We leave the precise statement for Section~$4$ (see Theorem
\ref{thm:moore2}). From the point of view of the theory of
integrable systems, the result is rather disappointing: there is
no such theory in connection with isothermic submanifolds of
dimension $n\geq 3$. The discrepancy between the richness of the
two-dimensional theory and the shortage of examples in higher
dimensions should be compared to that between the classes of
conformal maps between open subsets of Euclidean space in  two and
higher dimensions. In fact, our result contains as a special case
 Liouville's theorem on the classification of such conformal maps in  dimension $n\geq 3$.

 Having our definition of
isothermic submanifolds not led to as rich a theory as
    in the two-dimensional case, one might think that it should be more appropriate to choose a
    definition based on an extension of one of the properties that characterize isothermic
     surfaces in Euclidean three-space. For instance, one might look for those submanifolds
     that admit a nontrivial conformal sphere congruence, as suggested by Burstall, or
     search for the solutions of some natural extension  of  Christoffel's problem.
     We carry this out and give a complete classification of both classes of submanifolds.
     Unfortunately, they turn out to be distinct proper subclasses of the class of
     isothermic submanifolds as previously defined.

  Before we draw this introduction to a close, we provide  a brief guide to the contents of
  each section and describe some applications of our main result.

  We start Section~$2$ by reviewing from \cite{rs} the notions of nets and net morphisms
  between netted manifolds,
  a suitable setting for dealing with decomposition problems for manifolds and for maps
  between them.
  Then we recall  some  basic facts on twisted and warped products and state a characterization
  of twisted and warped product metrics on product manifolds in terms of  geometric properties
   of its product net
  \cite{mrs}.  We close the section by stating some results
  of \cite{to}, where a geometric characterization of higher dimensional  isothermic
  nets in
  Riemannian manifolds  was obtained, leading to a  conformal version of the local de Rham
  theorem.

       In Section $3$ we  give a  brief description of the model of Euclidean space as an
       umbilical hypersurface
       of the light cone in Minkowski space, and summarize some  facts that are needed in
       the sequel.
        In particular, we recall how to construct the
       aforementioned conformal diffeomorphism
      $\Theta\colon\,\Hy^m(-c)\times \Sf^{N-m}(c)\to \R^N$ and
      prove  a  formula relating the second fundamental forms of a conformal immersion
      into Euclidean space and its associated isometric immersion
      into  the light cone.

     Our main results are contained in Section $4$. We obtain a local classification of
     isothermic Euclidean submanifolds
     of dimension $n\geq 3$, or equivalently, we describe all conformal immersions of a
     Riemannian product  of
     dimension $n\geq 3$ into Euclidean space whose  second fundamental forms are adapted
     to the product net of the
     manifold. We follow with several applications. First, we specialize to surfaces in
     Euclidean space and generalize
     to arbitrary codimension a classical result of Bonnet on surfaces whose curvature lines
     have constant geodesic
     curvature (Corollary \ref{cor:surf}).  Then, we classify all conformal representations
     of Euclidean space of
     dimension $n\geq 3$ as a Riemannian product, that is, we determine all conformal
     local  diffeomorphisms of a Riemannian
     product of dimension $n\geq 3$ onto an open subset of Euclidean space
     (Corollary \ref{cor:moore2}). Besides
     containing Liouville's theorem as a special case, this also provides a geometric proof
     of the known
     classification of conformally flat Riemannian products \cite{la}.
     We follow with  a decomposition theorem for isometric immersions of a
     twisted product of two
     Riemannian manifolds into Euclidean space (Corollary \ref{cor:conf2}),
     which contains as a special case
     the main result of \cite{dft} on the classification of Euclidean submanifolds that
     carry a Dupin principal normal vector field with umbilical conullity. As a consequence,
     we provide an  alternate statement (Corollary \ref{cor:ciclides}) of
     the well-known classification of the cyclides of Dupin of arbitrary dimension \cite{ce1}.

   In Section $5$, after stating some basic facts on Codazzi tensors
   we formulate an extension of Christoffel problem for Euclidean submanifolds
   of arbitrary dimension and codimension, and determine all of its
   solutions of dimension $n\geq 3$.
   We conclude the paper by classifying in the last section all Euclidean submanifolds
   of dimension $n\geq 3$  that admit nontrivial conformal sphere congruences.

\section[Preliminaries]{Preliminaries}

    A suitable setting for treating decomposition results for manifolds
    was
    developed in \cite{rs} by introducing the category of   {\it netted
    manifolds\/}.   A {\it net\/} ${\cal E}=(E_i)_{i\in I_k}$ on a connected
    $C^{\infty}$-manifold $M$ is a splitting
    $TM=\oplus_{i\in I_k} E_i$ by a family
    of  integrable subbundles. Here and throughout the paper we denote $I_k=\{1,\ldots, k\}$.
    If $M$ is a Riemannian manifold and the subbundles are mutually orthogonal then the
    net is said to be an {\it orthogonal net\/}. The canonical net on a product manifold
$M=\Pi_{i=1}^k M_i$  is called the {\it product net\/}. A
$C^\infty$-map $\psi\colon\;M\to N$ between two {\it netted
manifolds\/}  $(M,{\cal E})$, $(N,{\cal F})$, that is, manifolds
$M,N$ equipped with nets ${\cal E}=(E_i)_{i\in I_k}$ and ${\cal
F}=(F _i)_{i\in I_k}$, respectively, is called a {\it net
morphism\/} if $d\psi (E_i(p))\subset F_i(\psi(p))$ for all $p\in
M$, $i\in I_k$, or equivalently, if for any $p\in M$ the
restriction $\psi|_{L_i^{{\cal E}}(p)}$ of $\psi$ to the leaf
$L_i^{{\cal E}}(p)$ of $E_i$ through $p$ is a $C^\infty$-map into
the leaf $L_i^{{\cal F}}(\psi(p))$ of $F_i$ through $\psi(p)$. The
net morphism $\psi$ is said to be a {\it net isomorphism\/} if in
addition it is a diffeomorphism and $\psi^{-1}$ is also a net
morphism. A net ${\cal E}$ on $M$ is said to be {\em locally
decomposable\/} if for every  $p\in M$ there exist a neighborhood
$U$ of $p$ in $M$ and a net isomorphism $\psi$ from $(U,{\cal
E}|_U)$ onto a product manifold $\Pi_{i=1}^k M_i$. The map
$\psi^{-1}\colon\;\Pi_{i=1}^k M_i\to U$ is called a {\em product
representation\/} of $(U,{\cal E}|_U)$.

     Given a product $M=\Pi_{i=1}^k M_i$ of $C^\infty$ manifolds $M_1,\ldots,M_k$,
      a metric $\<\;,\;\>$ on $M$ is called a {\em twisted product metric\/}
     if there exist Riemannian metrics
$\<\;,\;\>_i$ on $M_i$, $i\in I_k$, and a $C^\infty$ {\em
twist-function \/} $\rho=(\rho_1,\ldots,\rho_k)\colon\;M\to
\R^{k}_+$ such that
$\<\;,\;\>=\sum_{i=1}^k\rho_i^2\pi_i^*\<\;,\;\>_i.$ Then
$(M,\<\;,\;\>)$ is  said to be a {\em twisted product\/}  and is
denoted by $^\rho\Pi_{i=1}^k (M_i,\<\;,\;\>_i)$. When $\rho_1$ is
identically $1$ and $\rho_2,\ldots,\rho_k$ are independent of
$M_2,\ldots, M_k$, that is, there exist $\tilde{\rho}_i\in
C^\infty(M_1)$ such that $\rho_i=\tilde{\rho}_i\circ \pi_1$ for
$i=2,\ldots, k$, then $\<\;,\;\>$ is called a {\em warped product
metric\/} and
$(M,\<\;,\;\>):=(M_1,\<\;,\;\>_1)\times_{\tilde{\rho}}\Pi_{i=2}^k
(M_i,\<\;,\;\>_i)$ a {\em warped product\/} with {\em warping
function\/} $\tilde{\rho}=(\tilde{\rho}_2,\ldots,\tilde{\rho}_k)$.
If $\rho_i$ is identically $1$ for all $i\in I_k$, then
$\<\;,\;\>$ is a {\em Riemannian product metric\/} and
$(M,\<\;,\;\>)$  a {\em Riemannian product\/}.

An orthogonal net ${\cal E}=(E_i)_{i\in I_k}$ on a Riemannian
manifold  $M$ is called a $TP$-{\em net\/} if $E_i$ is umbilical
and $E_i^\perp$ is integrable for every $i\in I_k$. Recall that a
subbundle $E$ of $TM$ is {\em umbilical\/} if there exists a
vector field  $\eta$ in $E^\perp$ such that $(\nabla_X
Y)_{E^\perp}=\<X,Y\>\eta\;\;\mbox{for all }\;X,Y\in \Gamma(E).$
Here and in the sequel, the space of  smooth local sections of a
vector bundle $E$ over $M$ is denoted by $\Gamma(E)$, whereas
writing a vector field with a vector subbundle as a subscript
indicates taking the section of that vector subbundle obtained by
orthogonally projecting the vector field pointwise onto the
corresponding fiber of the subbundle. The vector field $\eta$ is
called the {\em mean curvature normal\/} of $E$. If, in addition,
$(\nabla_X \eta)_{E^\perp}=0$ for all $X\in \Gamma(E)$, then $E$
is said to be {\em spherical\/}. If $E$ is umbilical and its mean
curvature normal vanishes identically, then it is called {\it
totally geodesic\/} (or {\it auto-parallel\/}). An umbilical
distribution is automatically integrable, and the leaves are
umbilical submanifolds of $M$. When $E$ is totally geodesic or
spherical, its leaves are totally geodesic or spherical
submanifolds, respectively. By a {\it spherical submanifold\/} we
mean an umbilical \su whose mean curvature vector is parallel with
respect to the normal connection.

    An orthogonal net ${\cal E}=(E_i)_{i\in I_k}$ is called a $WP$-{\em net\/} if
    $E_i$ is spherical and
    $E_i^\perp$ is totally geodesic  for $i=2,\ldots, k$.
    This easily implies that  $E_1$ is  totally geodesic and
    $E_1^\perp$ is integrable, thus every $WP$-net is also a $TP$-net.
The terminologies $TP$-net and $WP$-net are justified by the
following result (see Proposition $4$ of \cite{mrs}).

\begin{proposition}\label{prop:tpwp}  On a connected product manifold
$M=:\Pi_{i=1}^k M_i$ the product net ${\cal E}=(E_i)_{i\in I_k}$
is a $TP$-net (resp., $WP$-net) with respect to a Riemannian
metric $\<\;,\;\>$ on $M$ if and only if $\<\;,\;\>$ is a twisted
product (resp., warped product) metric on $M$. Moreover, if
$\rho=(\rho_1,\ldots,\rho_k)\colon\;M\to \R^{k}_+$ is the twist
function and $U_i=-\grad\, (log\circ \rho_i)$, $i\in I_k$, where
the gradient is calculated with respect to $\<\;,\;\>$, then the
mean curvature normal of $E_i$ is $\eta_i=(U_i)_{E_i^{\perp}}$ for
all $i\in I_k$.
\end{proposition}

  An orthogonal net  ${\cal E}=(E_i)_{i\in I_k}$ on a Riemannian manifold is
  a {\em conformal product net},
  or a $CP$-net  for short, if
for $i=1,\ldots, k$ it holds that
  \be\label{eq:cpnets}E_i \mbox{ and } E_i^\perp \mbox{ are
umbilical and
}\<\nabla_{X_{\perp_i}}\eta_i,X_i\>=\<\nabla_{X_i}H_i,X_{\perp_i}\>\end{equation}
for all $X_i\in \Gamma(E_i)$ and  $X_{\perp_i}\in
\Gamma(E_i^\perp)$,
 where $H_i$ and $\eta_i$ are the mean curvature
normals of $E_i$ and $E_i^\perp$, respectively. We have (see
Proposition $4.2$ of
 \cite{to}):

\begin{proposition} \label{prop:conf}  On a connected
and simply connected product manifold $M=\Pi_{i=1}^k M_i$ the
product net ${\cal E}=(E_i)_{i\in I_k}$ is a $CP$-net with respect
to a Riemannian metric $\<\;\,,\;\>^\sim$ on $M$ if and only if
$\<\;\,,\;\>^\sim$ is conformal to a Riemannian product metric.
\end{proposition}

   By means of Proposition \ref{prop:conf}, the following conformal version of the local
   de Rham Theorem was obtained in \cite{to}.
   It shows that
   conformal product nets are natural generalizations of isothermic nets on surfaces.

\begin{theorem} \label{cor:drh} If a Riemannian manifold $M$ carries a $CP$-net
 ${\cal E}=(E_i)_{i\in I_k}$, then for every  $p\in M$ there exists a local
product representation $\psi\colon\;\Pi_{i=1}^k M_i\to U$ of
${\cal E}$ with $p\in U\subset M$, which is  conformal
 with respect to a  Riemannian product metric on
$\Pi_{i=1}^k M_i$.
\end{theorem}

\section[M\"obius geometry in the light cone]{M\"obius geometry in the light cone}

    We give a brief description of the model of Euclidean space
    as an umbilical hypersurface of the light cone of Minkowski
    space, a convenient setting for dealing with M\"obius
geometric notions. We refer the reader to \cite{hj} for  further
details.

    Let $\Les^{N+2}$ be the $(N\!+\!2)$--dimensional Minkowski space endowed with a
    Lorentz scalar product of signature $(+,\ldots,+, -)$, and let
$ \Ve^{N+1}= \{p\in\Les^{N+2}\colon\<p,p\>=0\} $ denote the light
cone in $\Les^{N+2}$. Then
$\Ee^N=\Ee^N_w=\{p\in\Ve^{N+1}\colon\<p,w\>=1\}$ is a model of
$N$--dimensional Euclidean space for any $w\in\Ve^{N+1}$. Namely,
choose $p_0\in \Ee^N$ and  a linear isometry $A\colon\,\R^N\to
\spa\{p_0,w\}^\perp$.
 Then the triple $(p_0,w,A)$ gives rise to an isometry
 $\Psi=\Psi_{p_0,w,A}\colon\;\R^N\to \Ee^N\subset \Les^{N+2}$ defined by
 $x\in \R^N\mapsto  p_0+A(x)-(1/2)\|x\|^2w$.

   Hyperspheres can be nicely described in $\Ee^N$: given a  hypersphere $S\subset \Ee^N$
   with (constant) mean curvature $h$ with respect to a  unit normal vector field
   $n$ along $S$, then $v=n_p+hp\in \Les^{N+2}$
is a constant space-like unit vector such that $\<v,p\>=0$ for all
$p\in S$; thus $S=\Ee^N\cap \{v\}^\perp$. Since $h=\<v,w\>$, then
$S$ is a hyperplane  if and only if $\<v,w\>=0$.

    The intersection angle of two (oriented) hyperspheres has also a simple description in
    this model: given hyperspheres $S_i=\Ee^N\cap \{v_i\}^\perp$ with unit normal vectors
    $n^i_p$, $1\leq i\leq 2$, at a common point $p$, their intersection angle at $p$ is
    given by $\<n^1_p,n^2_p\>=\<v_1,v_2\>$.  Thus $S_1$ and $S_2$ intersect orthogonally
    if and only if $\<v_1,v_2\>=0$.

    A hypersphere $S=\Ee^N\cap \{v\}^\perp$ has (Euclidean) center at $q_0\in \Ee^N$
    and mean curvature
    $h\neq 0$ if and only if $v=h q_0+(2h)^{-1}w$. This follows from $\<v,w\>=h$,
    $\<v,v\>=1$ and $\spa \{q_0,w\}^\perp \subset v^\perp$, the latter being due to the
    fact that any hyperplane through $q_0$ is orthogonal to $S$.

    Given a conformal immersion $G\colon\;M^n\to \Ve^{N+1}$
   with conformal factor $\va\in C^\infty(M^n)$, for any $\mu\in C^\infty(M^n)$ the map
   $G_\mu\colon\;M^n\to \Ve^{N+1}$,  $ p\mapsto \mu(p)G(p)$, is also conformal with conformal
   factor $\mu\va$.
   Therefore,  any conformal immersion $f\colon\;M^n\to \R^{N}$ with conformal factor
   $\va\in C^\infty(M^n)$ gives rise to an isometric immersion
   ${\cal I}(f)={\cal I}_{p_0,w,A}(f):=(\Psi\circ f)_{\va^{-1}}\colon\;M^n\to
   \Ve^{N+1}$. Conversely, if $F\colon\;M^n\to \Ve^{N+1}$ is an isometric immersion
   with $F(M^n)\subset \Ve^{N+1}\setminus\R_w$, where $\R_w=\{tw:t>0\}$,
   define ${\cal C}(F)={\cal C}_{p_0,w,A}(F)\colon M^n\to{\mathbb {R}}^{N}$   by
   $ \Psi\circ {\cal C}(F) = \Pi\circ F$, where
   $\Pi=\Pi_{w}\colon\,\Ve^{N+1}\setminus {\mathbb {R}}_w \to \Ee^{N}$
is the projection onto ${\mathbb {E}}^N$ given by
$\Pi(x)=x/\<x,w\>$. Since $\Pi$ is conformal with conformal factor
$\va_{\Pi}(x)=\<x,w\>^{-1}$, then ${\cal C}(F)$ is also conformal
with conformal factor $\va_{\Pi}\circ F=\<F,w\>^{-1}$.

In particular, conformal transformations of $\R^N$ are linearized
in this model: any $T\in \Oes_1(N+2)$ gives rise to a conformal
(M\"obius) transformation ${\cal T}={\cal C}(T\circ \Psi)$ of
$\R^{N}$ and, conversely, any M\"obius transformation of $\R^{N}$
is given in this way by means of some $T\in \Oes_1(N+2)$. For
instance, if $R$ is the reflection $R(p)=p-2\<p,v\>v$ with respect
to the hyperplane in $\Les^{N+2}$ orthogonal to the unit
space-like vector $v$, $\<v,w\>\neq 0$, then $I={\cal C}(R\circ
\Psi)$  is the inversion with respect to the hypersphere
$S=\Ee^N\cap \{v\}^\perp$. If $T\in \Oes_1(N+2)$ satisfies
$T(w)=\lambda w$ for some $\lambda\in \R$ then there exists a
similarity ${\cal H}$ of $\R^N$ of ratio $\lambda$ such that
${\cal C}(T\circ \Psi)={\cal H}$, i.e., $\Psi\circ {\cal H}
=\lambda  T\circ \Psi$. In particular, the isometries of $\R^N$
are given by those $T$  that fix $w$.

     Clearly, we have ${\cal C}_{p_0,w,A}({\cal I}_{p_0,w,A}(f))=f$ and ${\cal
I}_{p_0,w,A}({\cal C}_{p_0,w,A}(F))=F$  for any conformal
immersion $f\colon\,M^n\to \R^{N}$ and for any isometric immersion
$F\colon\;M^n\to \Ve^{N+1}$ with $F(M^n)\subset
\Ve^{N+1}\setminus\R_w$. For distinct triples $(p_0,w,A)$ and
$(\bar{p}_0,\bar{w},\bar{A})$,  there exist an inversion $I$ with
respect to a sphere of unit radius and a  similarity ${\cal H}$
such that \be\label{eq:comp}{\cal
C}_{p_0,w,A}(\Psi_{\bar{p}_0,\bar{w},\bar{A}}) =I\circ {\cal
H}.\ee Let us check this: consider the reflection
$R(p)=p-2\<p,v\>v$ determined by the unit space-like vector
$v=\<\bar{w},w\>^{-1}\bar{w}+(1/2)w$ and let $T\in\Oes_1(N+2)$ be
defined by $T(w)=R(\bar{w})=-(1/2)\<\bar{w},w\>w$,
$T(p_0)=R(\bar{p}_0)$ and $T\circ A=R\circ \bar{A}$. Then $R\circ
T$ takes $w$ to $\bar{w}$, $p_0$ to $\bar{p}_0$ and $R\circ T\circ
A= \bar{A}$, whence $ \Psi_{\bar{p}_0,\bar{w},\bar{A}}=R\circ
T\circ \Psi_{p_0,w,A}.$ Since ${\cal
C}_{p_0,w,A}(R\circ\Psi_{p_0,w,A})=I$ and ${\cal
C}_{p_0,w,A}(T\circ\Psi_{p_0,w,A})={\cal H}$  for an inversion $I$
in $\R^N$ with respect  to the hypersphere $S=\Ee_w^N\cap
\{v\}^\perp$ of unit radius  with center at
$\<\bar{w},w\>^{-1}\bar{w}\in \Ee_w^N$ and a  similarity ${\cal
H}$ on $\R^N$ of ratio $\lambda=-(1/2)\<\bar{w},w\>$, we obtain
$$\begin{array}{l}\Psi_{p_0,w,A}\circ {\cal
C}_{p_0,w,A}(\Psi_{\bar{p}_0,\bar{w},\bar{A}})=\Pi_{w}\circ
\Psi_{\bar{p}_0,\bar{w},\bar{A}}=\Pi_{w}\circ R\circ T\circ
\Psi_{p_0,w,A}\vspace{1ex}\\\hspace*{19.6ex}=\Pi_{w}\circ R\circ
\Psi_{p_0,w,A} \circ {\cal H}=\Psi_{p_0,w,A} \circ I\circ {\cal
H}.
\end{array}$$

Given a time-like vector $v\in \Les^{N+2}$ with $\<v,v\>=-1/c$ and
a linear isometry $B\colon\R^{N+1}\to \{v\}^\perp$, the isometric
immersion $T_{B,v}\colon\,\Sf^N(c)\to \Les^{N+2}$,
$X\in\Sf^N(c)\mapsto B(X)+v$, takes values in $\Ve^{N+1}$. When
$c=1$ and  $(p_0,w,A)$ is a triple such that $\{v,u\}$ is an
orthonormal basis of $\spa\{p_0,w\}$ with $w=v+u$ and
$A=B|_{\R^N}$, then a direct verification shows that ${\cal
C}_{p_0,w,A}(T_{B,v})$ is a stereographic projection of $\Sf^N$
onto $\R^N$. In general, there exists a similarity ${\cal H}$ of
$\R^{N+1}$ and a stereographic projection ${\cal
P}\colon\,\Sf^N\to \R^N$ such that \be\label{eq:sproj}{\cal
C}_{p_0,w,A}(T_{B,v})={\cal P}\circ {\cal H}.\ee

 Now we  construct a conformal
diffeomorphism
     $\Theta\colon\,\Hy^m(-c)\times \Sf^{N-m}(c)\to \R^N$ onto the
complement of an $(m-1)$-dimensional sphere  as  referred to in
the introduction. Namely, given an orthogonal decomposition
$\Les^{N+2}=V\oplus W$ with $V$ time-like and linear isometries
$C\colon\,\Les^{m+1}\to V$ and $D\,\colon\, \R^{N-m+1}\to W$,
define an isometric immersion $L_{C,D}\colon\,\Hy^m(-c)\times
\Sf^{N-m}(c)\to \Ve^{N+1}\subset \Les^{N+2}$ by  $(X,Y)\mapsto
C(X)+D(Y)$,
 and set \be\label{eq:Theta}\Theta={\cal C}_{p_0,w,A}(L_{C,D})\ee for some triple $(p_0,w,A)$.
 Observe that the
  image of $\Theta$ omits the
$(m-1)$-dimensional sphere  that is mapped onto $\Ee_w^N\cap V$ by
the isometry $\Psi_{p_0,w,A}\colon\;\R^N\to \Ee^N$. Alternately,
let $\Phi\colon\;\R_+^{m}\times_\sigma \Sf^{N-m}(c)\to
\R^N\setminus \R^{m-1}$ be the isometry given by $(X,Y)\mapsto
(x_1,\ldots,x_{m-1},\sigma(X)Y)$ for $X\in \R_+^{m}$ and $Y\in
\Sf^{N-m}(c)$, where
 $\R^m_+$ (resp., $\R^{m-1}$) denotes
the subset of points $X=(x_1,\ldots,x_m)$ of $\R^m$ where $x_m>0$
(resp., $x_m=0$), and $\sigma(X)=c^{1/2}x_m$. Endowing $\R_+^{m}$
with the metric $\sigma^{-2}\<\;,\;\>$, where $\<\;,\;\>$
  is the Euclidean metric, we obtain the half-space model of $\Hy^{m}(-c)$,
  and then $\Phi$ gives rise to a conformal diffeomorphism $\Theta_{\Phi}$ with conformal
  factor $\sigma\circ \pi_0$
  with respect to the product metric of $\Hy^m(-c) \times
  \Sf^{N-m}(c)$. Any $\Theta$ just defined differs from such a $\Theta_{\Phi}$ by an inversion in $\R^N$.

  We conclude this section by computing the relation
  between the second fundamental forms of a conformal immersion $f\colon\;M^n\to
  \R^{N}$ with conformal factor $\va\in C^\infty(M^n)$
   and its associated isometric immersion
  $F={\cal I}_{p_0,w,A}(f)\colon\;M^n\to
   \Ve^{N+1}\subset \Les^{N+2}$.

\begin{lemma} \label{le:sfs} The normal bundle $T_F^\perp M$ of
$F\colon\;M^n\to \Ve^{N+1}\subset \Les^{N+2}$ decomposes
orthogonally as $T_F^\perp M=d\Psi(T_f^\perp M)\oplus \Les^2$,
where $\Les^2$ is the Lorentzian plane bundle spanned by the
position vector $F$ and  $\eta:=\va w-d(\Psi\circ f)(\grad
\,\va^{-1})$, and the second fundamental form of $F$ splits
accordingly as \be\label{eq:sF}\alpha_F(X,Y)=\va\,\hess\,
\va^{-1}(X,Y)F+\va^{-1}d\Psi(\alpha_f(X,Y))-\<X,Y\>\eta,\ee where
$\hess$ and $\grad$ are calculated with respect to the metric
$\<\;\,,\;\>$ induced by $F$.
\end{lemma}
\proof Differentiating $F=\va^{-1}(\Psi\circ f)$ gives
$dF(X)=X(\va^{-1})(\Psi\circ f)+\va^{-1}d(\Psi\circ f)(X),$ which
already implies that $T_F^\perp M$ splits as stated. Now we
compute
\be\label{eq:1}\alpha_F(X,Y)=\bar{\nabla}_YdF(X)-dF(\nabla_Y
X),\ee where $\bar{\nabla}$ and $\nabla$ stand, respectively,  for
the (pull-back to $M^n$ of the) derivative in $\Les^{N+2}$ and the
Levi-Civita connection of $M^n$ with respect to $\<\;\,,\;\>$.
Using that
$$\bar{\nabla}_Y d(\Psi\circ
f)(X)=-\va^2\<X,Y\>w+d\Psi(\hat{\nabla}_Y df(X)),$$ where
$\hat{\nabla}$ is the (pull-back to $M^n$ of the) derivative on
$\R^N$, we have
\be\label{eq:2}\begin{array}{l}\bar{\nabla}_YdF(X)=YX(\va^{-1})(\Psi\circ
f)+X(\va^{-1})d(\Psi\circ
f)(Y)+Y(\va^{-1})d(\Psi\circ f)(X)\vspace{1ex}\\
\hspace*{13ex}-\va\<X,Y\>w+\va^{-1}d\Psi(\hat{\nabla}_Y df(X)).
\end{array}\ee
On the other hand, since
$$\nabla_Y X=\tilde{\nabla}_Y X+\va X(\va^{-1})Y+ \va
Y(\va^{-1})X-\va\<X,Y\>\,\grad \va^{-1},$$ where $\tilde{\nabla}$
is the Levi-Civita connection of the metric induced by $f$, we
obtain \be\label{eq:3}\begin{array}{l}dF(\nabla_Y X)=\nabla_Y
X(\va^{-1})(\Psi\circ f)+\va^{-1}d(\Psi\circ
f)(\nabla_Y X)\vspace{1ex}\nonumber\\
\hspace*{10.8ex} =\nabla_Y X(\va^{-1})(\Psi\circ
f)+\va^{-1}d\Psi(df(\tilde{\nabla}_Y
X))+X(\va^{-1})d(\Psi\circ f)(Y)\vspace{1ex}\\
\hspace*{13ex}+Y(\va^{-1})d(\Psi\circ f)(X)- \<X,Y\>d(\Psi\circ
f)(\,\grad \va^{-1}).\end{array}\ee Equation (\ref{eq:sF}) now
follows from (\ref{eq:1}), (\ref{eq:2}) and  (\ref{eq:3}).\qed

\section[Conformal immersions of Riemannian products]{Conformal immersions of
Riemannian products}

 In this section we prove our main result, namely,
   we give a complete description of all conformal immersions of a Riemannian product
   of dimension $n\geq 3$ into Euclidean space  whose second fundamental forms are
   adapted to the product net of the manifold.  As discussed in the introduction, this yields
   a local classification of all
   isothermic submanifolds of dimension $n\geq 3$ of Euclidean
   space.

  For an orthogonal decomposition $\R^N=\oplus_{i=1}^{k+1}\R^{m_i}$,
  with $\R^{m_{k+1}}$ possibly trivial, for every $i\in I_k$ let $N_i$
 denote either $\R^{m_i}$ or
 $\Sf^{m_i-1}(c_i)=\{X_i\in \R^{m_i} : \<X_i,X_i\>=1/c_i\}$ (in which case $m_i\geq 2$), and
 define $\Psi\colon\,\Pi_{i=1}^k N_i\to\R^N$ by $\Psi(x_1,\ldots,x_k)=
 (i_1(x_1),\ldots, i_k(x_k),v_{k+1})$, where
 $i_j\colon\,N_j\to \R^{m_j}$, $1\leq j\leq k$, is either the identity or the inclusion map,
 respectively, and $v_{k+1}\in \R^{k+1}$ is a constant vector. When  $N_i=\Sf^{m_i-1}(c_i)$ for
 every $i\in I_k$, then $\Psi$ takes values in (a small
 sphere of) $\Sf^{N-1}(c)$, with $1/c=\sum_{i=1}^k 1/c_i+\<v_{k+1},v_{k+1}\>.$ We call
 $\Psi$ an {\em extrinsic product\/} of $N_1,\ldots, N_k$.
 Given a product manifold $M=\Pi_{i=1}^k M_i$ we denote by
 $T_j\colon\;M\to M_j\times \Pi_{i\in I_k^j} M_i$ the map
 $(p_1,\ldots,p_k)\mapsto
 (p_j,(p_1,\ldots,\hat{p}_j,\ldots,p_k))$, where $I_k^j=I_k\setminus\{j\}$ and
the hat indicates that $p_j$ is missing.

\begin{theorem} \label{thm:moore2} Let $f\colon\;M^n:=\Pi_{i=1}^k M^{n_i}_i\to \R^{N}$,
$n\geq 3$, be a conformal immersion of a Riemannian product whose second fundamental form
is adapted to the product net of $M^n$.  Then one of the following  holds:
\begin{itemize}
\item[{\em (i)}] There exist an extrinsic product
$\Psi\colon\,\Pi_{i=1}^k N_i\to \Q_c^N\subset \R^{N+\epsilon}$ of
complete spherical  \sus of $\Q_c^N$, where $c\geq 0$,
$\epsilon=1$ if $c>0$ and  $\epsilon=0$ if $c=0$,  and \iis
$f_i\colon\;M_i\to N_i$, $1\leq i\leq k$, such that $$f={\cal
P}\circ H\circ\Psi\circ (f_1\times\cdots\times f_k),$$ where $H$
is a homothety of $\R^{N+\epsilon}$ and ${\cal P}$ is either a
stereographic projection of $H(\Q_c^N)=\Sf^N$ onto $\R^N$ if $c>0$
or an inversion in $\R^N$ with respect to a sphere of unit radius
otherwise. \vspace*{-1ex}

\begin{picture}(150,180)
\put(184,103){$\Q_c^N\subset \R^{N\!+\!\epsilon}$}
\put(31,127){$f_1$} \put(86,127){$f_k$}
\put(310,103){$H(\Q_c^N)\subset \R^{N\!+\!\epsilon}$}
\put(322,150){$\R^N$} \put(146,107){$\Psi$} \put(330,127){${\cal
P}$} \put(45,142){\vector(0,-1){25}}
\put(327,117){\vector(0,1){25}} \put(100,142){\vector(0,-1){25}}
\put(37,150){$M_1\!\times\cdots\times M_k$}
\put(37,103){$N_1\times\cdots\times N_k$}
\put(120,153){\vector(1,0){195}} \put(120,103){\vector(1,0){60}}
\put(244,103){\vector(1,0){60}}
\put(266,107){$H$}
\put(218,160){$f$}
\end{picture}
\vspace*{-20ex}

\item[{\em (ii)}]
For exactly one $j\in I_k$ there is an extrinsic product
$\Psi\colon\;\Pi_{i\in I_k^j}\, \Sf^{m_i}(c_i)\to\Sf^{N-m}(c)$,
\iis $f_j\colon\;M_j\to \Hy^{m}(-c)$  and $f_i\colon\;M_i\to
\Sf^{m_i}(c_i)$, $i\in I_k^j$, and a conformal diffeomorphism
$\Theta\colon\,\Hy^m(-c)\times \Sf^{N-m}(c)\to \R^N$    such that
$$ f =\Theta\circ(f_j\times
(\Psi\circ(f_0\times\cdots\times\hat{f_j}\times\cdots\times
f_k))\circ T_j).$$

\vspace*{-3ex}

\begin{picture}(150,184)
\put(10,110){Case $(ii)$ } \put(10,95){for $j=1$ }
\put(100,30){$\Hy^m(-c)\times \Sf^{N-m}(c)$} \put(95,95){$f_1$}
\put(121,127){$f_2$} \put(180,127){$f_k$} \put(315,149){$\R^N$}
\put(205,35){\vector(1,1){105}} \put(155,60){$\Psi$}
\put(269,76){$\Theta$} \put(108,142){\vector(0,-1){94}}
\put(135,142){\vector(0,-1){25}} \put(195,142){\vector(0,-1){25}}
\put(169,88){\vector(0,-1){40}} \put(97,150){$M_1\times
M_2\times\cdots\times M_k$}
\put(120,100){$\Sf_2(c_2)\!\times\!\cdots\!\times \!\Sf_k(c_k)$}
\put(205,153){\vector(1,0){100}} \put(250,160){$f$}
\end{picture}
\end{itemize}
\end{theorem}
\vspace*{-5ex} \proof  Define $F={\cal
I}_{p_0,w,A}(f)\colon\;M^n\to \Ve^{N+1}\subset \Les^{N+2}$ as in
Section $3$, and let ${\cal E}=(E_i)_{i\in I_k}$ be the product
net of $M^n$. We first prove:

\begin{lemma} \label{le:adap} The second fundamental form of $F$
is adapted to ${\cal E}$.
\end{lemma}
\proof It suffices to consider the case $k=2$ and then,
relabelling if necessary, we take $n_1\geq 2$.  Fixed
$p=(p_1,p_2)\in M^n$ and $\hat{X}\in E_2(p)$,  denote
$L=M_1^{n_1}\times \{p_2\}$, let $\bar{X}=d\pi_2(p)(\hat{X})\in
T_{p_2}M_2^{n_2}$ and, for any $q\in L$, let $\hat{X}(q)$ be the
unique vector in $E_2(q)$ that projects to $\bar{X}$ by
$d\pi_2(q)$. Then $\hat{X}$ is a parallel vector field along $L$
with respect to (the pull-back to $L$ of) the Levi-Civita
connection of $M^n$. Let $\xi=dF(\hat{X})\colon\,L\to \Les^{N+2}$.
Then for any  $X\in TL$ we have  \be \label{eq:dxi}
d\xi(X)=\alpha_F(X,\hat{X})=\omega(X)F,\,\,\,\,\mbox{with}\,\,\omega(X)=\va\,\hess\,
\va^{-1}(X,\hat{X}),\ee where the second equality follows from
(\ref{eq:sF}) and the assumption that  the second fundamental form
of $f$ is adapted to ${\cal E}$.
 For $X,Y\in T_pL$ linearly independent,
the exterior derivative of (\ref{eq:dxi}) gives
$$0=d^2\xi(X,Y)=d\omega(X,Y)F-\omega(X)dF(Y)-\omega(Y)dF(X).$$
Since  $F$, $dF(X)$ and $dF(Y)$ are linearly independent, because
 $F$ is an immersion and the position vector $F$ is a nonzero
normal vector field, we conclude that $\omega$ vanishes  whence
$\alpha_F(X,\hat{X})=0$ for any $X\in TL$. \vspace{1ex}\qed

\begin{lemma} \label{le:ortog} Let $p, q\in M^n$, $X_i\in E_i(p)$ and $X_j\in E_j(q)$,
$i\neq j$. Then $dF(X_i) \perp dF(X_j)$.
\end{lemma}
\proof Let $\bar{X}_i=d\pi_i(p)(X_i)\in T_{\pi_i(p)}M_i^{n_i}$
and, for  any point $z$ of the fiber $M_{\perp_i}$ of $\pi_i$
through $p$, let $X_i(z)$ be the unique vector in $E_i(z)$ that
projects to $\bar{X}_i$ by $d\pi_i(z)$. Arguing as in the
beginning of the proof of Lemma \ref{le:adap} and using its
conclusion, we obtain that  $dF(X_i)\colon\,M_{\perp_i}\to
\Les^{N+2}$ is constant, and the statement
follows.\vspace{1ex}\qed

  Now, for $i\in I_k$ define linear subspaces $W_i$ of
  $\Les^{N+2}$ by
  $$W_i=\spa \{dF_q(X_i) : q\in M^n, X_i\in E_i(q)\}.$$
By Lemma \ref{le:ortog} the subspaces $W_i$ are mutually
orthogonal. We distinguish two cases.\vspace{1ex}

First suppose that all the $W_i$ inherit non-degenerate metrics
from $\Les^{N+2}$. Define $W_{k+1}=(W_1\oplus\cdots \oplus
W_k)^\perp$ and let $P_i\colon\,\Les^{N+2}\to W_i$ denote
orthogonal projection. Since  $W_1, \ldots, W_{k+1}$ are mutually
orthogonal, for any $X_j\in \Gamma(E_j)$ we have $d(P_i\circ
F)(X_j)=0$ for $i\neq j$. Thus,  for $i\leq k$, $P_i\circ F$ is
constant on the fibers of the projection $\pi_i\colon\,M^n\to
M_i^{n_i}$, while $P_{k+1}\circ F$ has a constant value $e_{k+1}$
on $M^n$. Fixed $\bar{p}=(\bar{p}_1,\ldots,\bar{p}_{k})\in M^n$,
for $i\in I_k$ define $F_i\colon\,M_i^{n_i}\to W_i$ by
$F_i=P_i\circ F\circ i_{\bar{p}}$, where
$i_{\bar{p}}\colon\,M_i^{n_i}\to M^n$ denotes the isometric
inclusion of $M_i^{n_i}$ into $M^n$ given by  $p_i\mapsto
(\bar{p}_1,\ldots,p_{i},\ldots,\bar{p}_{k})$.  Then $F_i$ is an
isometric immersion with $P_i\circ F=F_i\circ \pi_i$ for every
$i\in I_k$, whence\be \label{eq:F} F=\sum_{i=1}^{k+1}P_i\circ
F=\sum_{i=1}^k F_i\circ \pi_i +e_{k+1}.\ee From (\ref{eq:F}) we
get
$$0=\<F,F\>=\sum_{i=1}^k\<F_i,F_i\>\circ
\pi_i+\<e_{k+1},e_{k+1}\>,$$ from which we conclude that each
$\<F_i,F_i\>$ is a constant, say, $1/c_i$. Notice that $W_j$ is
time-like for exactly one $j\in I_{k+1}$. Assume first that
$j=k+1$. Orthogonally decompose $e_{k+1}=v+u$ with $v\in W_{k+1}$
time-like, set $-1/c=\<v,v\>$ and
$\hat{F}=\sum_{i=1}^k\<F_i,F_i\>\circ \pi_i+u$, whence
$\<\hat{F},\hat{F}\>=1/c$.  Then there exist an extrinsic product
$\Psi\colon\;\Pi_{i=1}^k \Sf^{m_i-1}(c_i)\to\Sf^{N}(c)\subset
\R^{N+1}=\oplus_{i=1}^{k+1}\R^{m_i}$,  linear isometries
$B_i\colon\,\R^{m_i}\to W_i$, $i\in I_k$,
$B_{k+1}\colon\,\R^{m_{k+1}}\to \hat{W}_{k+1}:= W_{k+1}\cap
v^\perp$ and $B=\oplus_{i=1}^{k+1} B_i
\colon\,\R^{N+1}=\oplus_{i=1}^{k+1}\R^{N_i}\to \oplus_{i=1}^{k}
W_i\oplus \hat{W}_{k+1}=v^\perp$  such that $\hat{F}=B\circ
\Psi\circ (f_1\times \cdots\times f_k)$, where
$f_i\colon\,M_i^{n_i}\to \R^{m_i}$ is defined by $F_i=B_i\circ
f_i$ for every $i\in I_k$. Therefore $F=T_{B,v}\circ \Psi\circ
(f_1\times \cdots\times f_k)$, and we obtain using
(\ref{eq:sproj}) that
$$\begin{array}{l}f={\cal C}_{p_0,w,A}({\cal I}_{p_0,w,A}(f))={\cal
C}_{p_0,w,A}(F)={\cal C}_{p_0,w,A}(T_{B,v})\circ \Psi\circ
(f_1\times \cdots\times f_k)\vspace{1ex}\\\hspace{2ex}={\cal
P}\circ {\cal H}\circ\Psi\circ (f_1\times\cdots\times
f_k),\end{array}$$ where ${\cal H}$ is a similarity of $\R^{N+1}$
and ${\cal P}$ is  a stereographic projection of ${\cal
H}(\Sf^{N}(c))=\Sf^N$ onto $\R^N$. Writing ${\cal H}=H\circ S$ for
a homothety $H$ and an isometry $S$ of $\R^{N+1}$, and observing
that $S\circ \Psi$ is still an extrinsic product of $\Pi_{i=1}^k
\Sf^{m_i-1}(c_i)$ into $\Sf^{N}(c)$, we obtain case $(i)$ of the
statement for $c>0$.

 Now suppose that $j<k+1$. Choose linear isometries $C\colon\,\Les^{m+1}\to
 W_j$, $D_i\colon\,\R^{m_i}\to W_i$, $i\in I^j_{k+1}$, define
$f_j\colon\,M_j^{n_j}\to \Les^{m+1}$ by $F_j=C\circ f_j$ and
$f_i\colon\,M_i^{n_i}\to \R^{m_i}$ by $F_i=D_i\circ f_i$. Set
 $D=\oplus_{i\in I^j_{k+1}} D_i \colon\,\R^{N-m+1}=
 \oplus_{i\in I^j_{k+1}}\R^{m_i}\to \oplus_{i\in I^j_{k+1}} W_i$. Then there
exist an  extrinsic product $\Psi\colon\;\Pi_{i\in I^j_{k}}
\Sf^{m_i}(c_i)\to\Sf^{N-m}(c)$, $1/c=\sum_{i\in I^j_{k}}
1/c_i+\<e_{k+1},e_{k+1}\>,$ such that $F=L_{C,D}\circ \bar{f}$,
where $\bar{f}=(f_j\times
(\Psi\circ(f_1\times\cdots\times\hat{f_j}\times\cdots\times
f_k))\circ T_j).$ Therefore, defining a conformal diffeomorphism
$\Theta\colon\,\Hy^m(-c)\times \Sf^{N-m}(c)\to \R^N$ by
(\ref{eq:Theta}), we obtain
$$f={\cal
C}_{p_0,w,A}(F)={\cal C}_{p_0,w,A}(L_{C,D})\circ
\bar{f}=\Theta\circ\bar{f},$$ and we are in the situation of case
$(ii)$ of the statement.

   The second case we must consider is when some of the $W_i$
   inherit a degenerate metric from $\Les^{N+2}$. In this case,
   without loss of generality, we may assume that $W_1,\ldots,
   W_\ell$ have such a metric while $W_{\ell+1},\ldots, W_k$ have
   non-degenerate (and so necessarily space-like) metric. Then there exists a
   $1$-dimensional and light-like subspace $L_0$ such that
   $W_i\cap W_i^\perp=L_0$ for $i=1,\ldots, \ell$. Choose a
   second, distinct light-like line $L_1$ orthogonal to $W_{\ell+1},\ldots,
   W_k$. Set $\hat{W}_i=W_i\cap L_1^\perp$ (so that
   $\hat{W}_i=W_i$ for $i>\ell$) and finally set
   $$\hat{W}_{k+1}=(L_0\oplus \hat{W}_1\oplus\cdots \oplus
   \hat{W}_k\oplus L_1)^\perp.$$ We therefore have a decomposition
   $\Les^{N+2}=L_0\oplus \hat{W}_1\oplus\cdots \oplus
   \hat{W}_{k+1}\oplus L_1$
   and corresponding orthogonal projections
   $\hat{P}_i\colon\,\Les^{N+2}\to \hat{W}_i$. Arguing as in the preceding case, we obtain
   that, for
   $i\in I_k$,
   $\hat{P}_i\circ F$ is constant on the fibers of
   $\pi_i$   while the components of $F$ in $L_1$ and
   $\hat{W}_{k+1}$ are constant. Thus, there exist isometric
   immersions $\hat{F}_i\colon\,M_i^{n_i}\to \hat{W}_i$ such that
   \be\label{eq:F2}F=\bar{p}_0+\sum_{i=0}^k\hat{F}_i\circ \pi_i+e_{k+1}+\<F,\bar{p}_0\>\bar{w},\ee
   where  $\bar{p}_0$ is a
   light-like constant vector in $L_1$, $e_{k+1}$ is a (space-like) constant vector in
   $\hat{W}_{k+1}$ and $\bar{w}\in L_0$ is chosen so that $\<\bar{p}_0,\bar{w}\>=1$.
   From (\ref{eq:F2}) and
   $\<F,F\>=0$ we conclude that
   $2\<F,\bar{p}_0\>=-\sum_{i=1}^k\<\hat{F}_i\circ \pi_i, \hat{F}_i\circ
   \pi_i\>$. Identifying $(L_0\oplus L_1)^\perp$ with $\R^N$ by means of a linear
   isometry $\bar{A}$,
    the isometric immersion
   $\hat{F}=\sum_{i=1}^k\hat{F}_i\circ \pi_i+e_{k+1}$ is an extrinsic product of
   $\hat{F}_1,\ldots, \hat{F}_k$ with respect to the orthogonal decomposition
   $\R^N=\oplus_{i=1}^{k+1} \hat{W}_i$, whose image lies in the affine subspace
   $e_{k+1}\oplus \hat{W}^\perp_{k+1}$, and $F=
   \Psi_{\bar{p}_0,\bar{w},\bar{A}}\circ \hat{F}$. Therefore, using (\ref{eq:comp}) we conclude
   that there exist
   an inversion $I$ with respect to a sphere of unit
radius and a similarity ${\cal H}$ on $\R^N$ such that
   $$f={\cal C}_{p_0,w,A}(F)={\cal
   C}_{p_0,w,A}(\Psi_{\bar{p}_0,\bar{w},\bar{A}}\circ \hat{F})={\cal
   C}_{p_0,w,A}(\Psi_{\bar{p}_0,\bar{w},\bar{A}})\circ \hat{F}=I\circ {\cal H}\circ
   \hat{F},\,\,\,$$
which gives case $(i)$ of the statement for $c=0$. \qed

\begin{remarks}\po\label{re:main}{\em $(a)$ An alternate
statement for part $(ii)$
is as follows:\\
 $(ii')$
 For exactly one $j\in I_k$
 there exist an isometry
 $\Phi\colon\,\R_+^m\times_{\sigma}\Sf^{N-m}(c)\to \R^N$, an extrinsic product
  $\Psi\colon\;\Pi_{i\in I_k^j} \Sf^{m_i}(c_i)\to\Sf^{N-m}(c)$, a conformal
 immersion $f_j\colon\;M_j\to \R_+^m$ with conformal factor $(\sigma \circ f_j)$,
 \iis $f_i\colon\;M_i\to \Sf^{m_i}(c_i)$, $i\in I_k^j$,  and an inversion $I$  in $\R^N$
  such that $ f =I\circ \Phi\circ(f_j\times
(\Psi\circ(f_1\times\cdots\times\hat{f_j}\times\cdots\times
f_k)))\circ T_j.$\vspace{1ex}\\ $(b)$  Case $(i)$ with $c>0$ can
occur only if $k\leq N\!-\!n$. This follows from the fact the
codimension $N\!-\!n$ of $f$ is greater than or equal to the
codimension $k$ of $\Psi$. Similarly, in  case $(ii)$ we must have
$k\leq N\!-\!n\!+1$.}
\end{remarks}

    Clearly, Theorem \ref{thm:moore2} does not hold for $n=2$. In fact, the proof of
    Theorem \ref{thm:moore2} yields a classification of a special class of isothermic
    surfaces in $\R^N$, which extends  a theorem of
   Bonnet for $N=3$ (\cite{bo}, cf.  \cite{da2}, vol. III, p. 121)).

\begin{corollary} \label{cor:surf} Let $f\colon\;M^2\to  \R^N$ be a surface with \fnb
without umbilic points. Assume that the curvature lines of both families have constant
geodesic curvature. Then one of the following  holds:
\begin{itemize}
\item[{\em (i)}] There exist local isothermic parameterizations by curvature lines
$\psi\colon\; J_1\times J_2 \to M^2$, an isometry $\Phi\colon\;\R^{N_1}\times \R^{N_2}\to \R^N$,
 an inversion $I$ in $\R^N$ and unit speed curves $\alpha\colon\;J_1\to \R^{N_1}$ and
  $\beta\colon\;J_2\to \R^{N_2}$  such that $f\circ \psi=I\circ \Phi\circ (\alpha\times \beta)$.
\item[{\em (ii)}] There exist local isothermic parameterizations by curvature lines
$\psi\colon\; J_1\times J_2 \to M^2$, unit speed curves
$\alpha\colon\;J_1\to \Hy^{m}(-c)$ and $\beta\colon\;J_2\to
\Sf^{N-m}(c)$, and a conformal diffeomorphism
$\Theta\colon\,\Hy^m(-c)\times \Sf^{N-m}(c)\to \R^N$ such that $
f\circ \psi = \Theta\circ(\alpha\times \beta)$.
\end{itemize}
\end{corollary}
\proof Let ${\cal E}=(E_1,E_2)$ be the orthogonal net on $M^2$
determined by the principal directions of $f$. By  assumption,
$E_1$ and $E_2$ are spherical whence ${\cal E}$ is a $CP$-net  by
(\ref{eq:cpnets}). By Theorem \ref{cor:drh}, there exist local
isothermic parameterizations $\psi\colon\; J_1\times J_2 \to M^2$
whose coordinate curves are integral curves of $E_1$ and $E_2$.
Moreover, if $\va$ denotes the conformal factor of $\psi$, then it
is easily checked that the fact that  the coordinate curves of
$\psi$ have constant geodesic curvature is equivalent to $\hess
\va^{-1}$ (the Hessian being computed with respect to the flat
metric on $J_1\times J_2$) being adapted to the product net of
$J_1\times J_2$, i.e., $\va^{-1}=\va_1\circ \pi_1+ \va_2\circ
\pi_2$ for some $\va_1\in C^\infty(J_1)$ and $\va_2\in
C^\infty(J_2)$. Therefore Lemma \ref{le:adap}, whence the
remaining of the proof of Theorem \ref{thm:moore2}, applies for
$f\circ \psi$\vspace{1ex}\qed

      In the case $n=N$, Theorem \ref{thm:moore2}  gives a classification of
    all {\em conformal representations\/} of Euclidean space as a Riemannian product,
    which contains as a
    special case Liouville's theorem on the classification of conformal mappings between open
    subsets
    of $\R^N$, $N\geq 3$. It also provides a geometric proof of the classification of
    conformally flat Riemannian products
    (cf. \cite{la}, Section D, Proposition~2).

\begin{corollary} \label{cor:moore2} Let $f\colon\;M^N:=\Pi_{i=1}^k M^{N_i}_i\to \R^{N}$,
$N\geq 3$,  be a conformal local diffeomorphism of a Riemannian
product. Then one of the following holds:
\begin{itemize}
\item[{\em (i)}] There exist an isometry $\Phi\colon\;\Pi_{i=1}^k\R^{N_i}\to \R^N$,
local isometries $f_i\colon\;M^{N_i}_i\to \R^{N_i}$, an inversion
$I$ and a homothety $H$ in $\R^{N}$ such that $f=I\circ H\circ
\Phi\circ (f_1\times\cdots\times f_k).$
\item[{\em (ii)}] $k=2$ and, after relabelling the factors if necessary,  there exist local
isometries
$f_1\colon\;M^{N_1}_1\to \Hy^{N_1}(-c)$ and
$f_2\colon\;M^{N_2}_2\to \Sf^{N_2}(c)$ and a  conformal
diffeomorphism $\Theta\colon\,\Hy^{N_1}(-c)\times \Sf^{N_2}(c)\to
\R^N$ such that $ f =  \Theta\circ(f_1\times f_2).$
\end{itemize}
\end{corollary}

Theorem \ref{thm:moore2}  has also the following  consequence for
isometric immersions of twisted products into Euclidean space.

\begin{corollary} \label{cor:conf2} Let $f\colon\;(M^n,\<\;,\;\>):=M^{n_1}\times M^{n_2}\to
\R^N$, $n\geq 2$, be an \ii of a simply connected twisted product
whose second fundamental form  is adapted to the product net
${\cal E}=(E_1,E_2)$ of $M^n$. If $n_i=1$ for some $i\in \{1,2\}$,
suppose further that the leaves of $E_i$ are spherical in $M^n$.
Then $\<\;,\;\>$ is conformal to a product metric $\<\;,\;\>^\sim$
and the conclusion of Theorem \ref{thm:moore2} holds for
$f\colon\;(M^n,\<\;,\;\>^\sim)\to  \R^N$.
\end{corollary}
\proof For $n=2$ the result reduces to Corollary \ref{cor:surf},
thus we may assume that $n\geq 3$. By Proposition \ref{prop:tpwp},
 ${\cal E}$ is a $TP$-net. Then,
Lemma \ref{prop:ext2} below and the assumption for the case when
$n_i=1$ imply that $E_1$ and $E_2$ are spherical subbundles,
whence ${\cal E}$ is a $CP$-net by (\ref{eq:cpnets}). Thus
$\<\;,\;\>$ is conformal to a product metric $\<\;,\;\>^\sim$ by
Proposition \ref{prop:conf}. Since the second fundamental form of
$f$ is adapted to~${\cal E}$, Theorem \ref{thm:moore2} applies for
$f\colon\;(M^n,\<\;,\;\>^\sim)\to  \R^N$.

   \begin{lemma}\po \label{prop:ext2} Let $h\colon\,M^n\to\R^{N}$ and
$g\colon\,L^k\to M^{n}$, $k\geq 2$, be isometric immersions. If
$g$ is umbilical and $\alpha_h(dg(X),Z)=0$ for all $X\in
\Gamma(TL)$,   $Z\in \Gamma(T_g^\perp L)$,  then $g$ is spherical.
\end{lemma}
\proof  Since $g$ is umbilical,  the Codazzi equation for $g$
yields \be \label{eq:codg} \<Y,T\>\nabla^\perp_X
H_g-\<X,T\>\nabla^\perp_Y H_g=
 (R(dg(X),dg(Y))dg(T))_{T_g^\perp L} \ee  for
all $X,Y,T\in \Gamma(TL)$, where $R$ stands for  the curvature
tensor of $M^n$ and $H_g$ for the  mean curvature vector of $g$.
The  Gauss equation of $h$ and the assumption on $\alpha_h$ imply
that $R(dg(X),dg(Y))Z=A^h_{\alpha_h(dg(Y),
Z)}dg(X)-A^h_{\alpha_h(dg(X), Z)}dg(Y)=0$ for all $X,Y\in
\Gamma(TL)$ and $Z\in \Gamma(T_g^\perp L)$, whence the
right-hand-side of (\ref{eq:codg}) vanishes. Choosing $Y=T$
orthogonal to $X$ we conclude that $H_g$ is parallel in the normal
connection. \vspace{1ex} \qed

 An interesting particular case is that of a cyclide of Dupin $f\colon\,M^{n}\to\R^{n+1}$,
 in which case the
 eigenbundles correspondent to the two distinct principal curvatures give rise to a $TP$-net on
 $M^n$. Recall that $f$  is a {\em cyclide of Dupin\/}
 of characteristic $(m,n-m)$ if it has everywhere two distinct principal  curvatures of
multiplicities $m$ and $n-m$, respectively, which are constant
along the corresponding eigenbundles. One can now easily derive
the following alternate  classification of the cyclides of Dupin
(cf. \cite{ce1}).

\begin{corollary}\po\label{cor:ciclides}  Let $f\colon\,M^{n}\to\R^{n+1}$ be a cyclide of Dupin
 of characteristic
 $(m,n-m)$. Then there exist a conformal diffeomorphism $\psi$ of  an open
subset $W\subset \Q^{n-m}_{c}\times \Sf^{m}$, $c>-1$, onto $M^n$,
and a conformal diffeomorphism $\Theta\colon\,\Hy^{n-m+1}\times
\Sf^{m}\to \R^{n+1}$  such that $f\circ \psi= \Theta\circ
(f_1\times i)|_W$, where $f_1\colon\; \Q^{n-m}_{c}\to\Hy^{n-m+1}$
is a spherical inclusion and $i\colon\;\Sf^{m}\to \Sf^{m}$ is the
identity map. Moreover, the classes of conformally congruent
cyclides of Dupin are parameterized by their characteristic and
the value of $c$.
\end{corollary}
\begin{remarks}\po\label{re:cyclide}{\em $(i)$ If $c>0$ in Corollary \ref{cor:ciclides} then we
have an alternate description:\vspace{.5ex} \\ There exist a
conformal diffeomorphism $\psi$ of an open subset $W\subset
\Sf^{n-m}(c)\times \Sf^{m}$ onto $M^n$, an isometric embedding
$\Psi$ of $\Sf^{n-m}(c)\times \Sf^{m}$ into
$\Sf^{n+1}(\tilde{c})$, $\tilde{c}=c/(c+1)$, as an extrinsic
product, a homothety $H$ in $\R^{n+2}$ and a stereographic
projection $\P$ of $\Sf^{n+1}(\bar{c})=H(\Sf^{n+1}(\tilde{c}))$
onto $\R^{n+1}$ such
that $f\circ \psi=\P\circ H\circ \Psi|_{W}$.\vspace{1ex}\\
$(ii)$ Another particular case of Corollary \ref{cor:conf2} is
 the main result of \cite{dft}, which classifies Euclidean submanifolds that
 carry a Dupin principal normal  with umbilical conullity. }
\end{remarks}

\section[Christoffel problem
]{Christoffel problem }

    In this section we classify submanifolds that are solutions of a generalized
    Christoffel problem. We start by summarizing some facts on Codazzi tensors that will be
    needed in the sequel. Recall that a symmetric tensor ${\cal
    S}$ on a
Riemannian manifold $M$
      is  a {\it Codazzi tensor\/}  if
      $$(\nabla_X{\cal S})Y=(\nabla_Y{\cal S}) X\;\;\mbox{for all}\;\;
X,Y\in \Gamma(TM),$$ where $(\nabla_X{\cal S})Y=\nabla_X{\cal S}
Y-{\cal S}(\nabla_X Y)$. The following basic result is due to
Reckziegel \cite{re} (cf. Proposition $5.1$ of \cite{to}).

\begin{proposition} \label{prop:codtensor} Let ${\cal S}$ be a Codazzi tensor on a
Riemannian manifold $M$ and let $\lambda\in C^\infty(M)$ be an
eigenvalue of ${\cal S}$ such that $E_\lambda=\ker(\lambda I-{\cal
S})$ has constant rank $k$. Then:
\begin{itemize}
\item[{\em (i)}] If $k\geq 2$ then $\lambda$ is constant along $E_\lambda$.
\item[{\em (ii)}] $E_\lambda$ is an umbilical distribution (in fact spherical if $\lambda$ is constant
along $E_\lambda$) with mean curvature normal $\eta$ given by
$(\lambda I-{\cal S})\,\eta=(\nabla \,\lambda)_{E_\lambda^\perp}.$
\end{itemize}
\end{proposition}

 We now consider Codazzi tensors with exactly two distinct eigenvalues
  everywhere. The next result is contained in Theorem $5.2$ of \cite{to}.

\begin{proposition} \label{thm:codtensor} Let $M$ be a Riemannian manifold and let
${\cal S}$ be a Codazzi tensor on $M$ with exactly two distinct
eigenvalues $\lambda$ and $\mu$ everywhere. Let $E_\lambda$ and
$E_\mu$ be the corresponding eigenbundles. Then
$(E_\lambda,E_\mu)$ is a $CP$-net if and only if
$$ 2r X(h)Y(h)+h X(h)Y(r)+h Y(h)X(r)-r h\hess h(X,Y)=0, $$ where
$\lambda=h(1-r)$ and $\mu=h(1+r)$.
\end{proposition}

 The following is a special case of Corollary $5.3$ of \cite{to}. For a product manifold
 $M=M_1\times M_2$  with product net $(E_1, E_2)$, we denote by $\Pi_i\colon\;TM\to E_i$
 the canonical projection onto $E_i$, $1\leq i\leq 2$.

\begin{proposition}\po\label{prop:codazzi23} Let ${\cal S}$ be a Codazzi tensor on
 a Riemannian manifold $M$  with eigenvalues $\lambda$ and $-\lambda$, where
$\lambda\neq 0$ everywhere. Let $E_+$ and $E_-$ be the
corresponding eigenbundles and assume that $\lambda$ is constant
along $E_-$. Then one of the following  holds:
\begin{itemize}
\item[{\em (i)}] $\lambda$ is constant along $E_+$ and for every point
$p\in M$ there exists a local product representation
$\psi\colon\;M_1\times M_2\to U$ of $(E_+, E_-)$ with $p\in
U\subset M$, which is an isometry with respect to a Riemannian
product metric $\<\;,\;\>$ on $M_1\times M_2$. Moreover,
$(d\psi)^{-1}\circ {\cal S}\circ d\psi=a(\Pi_2-\Pi_1)$ for some
$a\neq 0$.
\item[{\em (ii)}] for every point $p\in M$ there exists a local product representation
$\psi\colon\;I\times N\to U$ of $(E_+, E_-)$ with $p\in U\subset
M$, where $I\subset \R$ is an open interval, which is an isometry
with respect to a warped  product metric $\<\;,\;\>$ on $I\times
N$ with warping function $\rho$.  Moreover, $(d\psi)^{-1}\circ
{\cal S}\circ d\psi=a(\rho\circ \pi_1)^{-2}(\Pi_2-\Pi_1)$ for some
$a\neq 0$.
\end{itemize}
\end{proposition}

    The classical proof that isothermic surfaces are precisely the ones that admit a
    {\em dual\/} surface (or  Christoffel transform) starts by showing that if two surfaces
    are mapped conformally onto each other with parallel tangent planes at corresponding
    points, then either they are a pair of minimal surfaces with the same conformal structures
     and orientations, or the correspondence between them necessarily preserves curvature lines.
     Then,  the proof proceeds by showing that in the latter case either the surfaces differ
     by a homothety and a translation or they are a Christoffel pair of isothermic surfaces.

   Classically, two surfaces $f\colon\;M^2\to \R^{3}$ and $\Fes\colon\,M^2\to\R^{3}$ that can
   be mapped onto each other with preservation of curvature lines and with parallel tangent
   planes at corresponding points are said to be related by a {\it
Combescure transformation}  (cf. \cite{bi}, v.II-1, p. 108).
Therefore, Christoffel's characterization of isothermic surfaces
implies that a surface is isothermic if and only if it admits a
nontrivial conformal Combescure transform. By nontrivial we mean
that the surfaces do not differ by a composition of a homothety
and a translation.

  The Combescure transformation can be extended  for
Euclidean submanifolds of arbitrary dimension and codimension as
follows (cf. \cite{dt2}).  Given an \ii $f\colon\;M^n\to \R^{N}$,
a map $\Fes\colon\,M^n\to\R^{N}$ is said to be
  a {\it Combescure transform\/} of $f$ determined by a symmetric tensor ${\cal S}$ on $M^n$
  if $d\Fes=df\circ{\cal S}$.
This implies (see Proposition~$1$ of \cite{dt2}) that  ${\cal S}$
is a {\em commuting\/} Codazzi tensor, i.e.,
$$ \a_f(X,{\cal S} Y)=\a_f({\cal S} X,Y)\fall X,Y\in \Gamma(TM^n).
$$ If ${\cal S}$ is invertible, then $\Fes$ is an immersion with
the same Gauss map as $f$ into the Grassmann manifold of non
oriented $n$--planes in $\R^{N}$. Moreover, the requirement that
the tensor ${\cal S}$ be symmetric reduces in the surface case to
the assumption that $f$ and $\Fes$ have the same curvature lines,
since their \sffs  are related by \be\label{eq:rsff}
\a_{\Fes}(X,Y)=\a_f({\cal S} X,Y) \fall X,Y\in \Gamma(TM^n). \ee
If $M^n$ is simply-connected and $\Fes\colon\,M^n\to\R^{N}$ is any
Combescure transform  of $f$ determined by a symmetric tensor
${\cal S}$, it was shown in Proposition $3$ of \cite{dt2} that
there exist $\va\in C^{\infty}(M^n)$ and $\beta\in
\Gamma(T_f^{\perp}M^n)$ satisfying \be\label{eq:gnorm} \a_f(\nabla
\va,X)+\nap_X \beta=0 \fall X\in \Gamma(TM^n), \ee such that
${\cal S}$ and $\Fes$ are given by \be\label{eq:fes} {\cal
S}={\cal S}_{\va, \beta}:=\hess\varphi - A^f_\beta
\;\;\;\;\mbox{and}\;\;\;\; \Fes=\Fes_{\va, \beta}:=df(\grad\,\va)
+\beta. \ee Conversely, if $\va\in C^{\infty}(M^n)$ and $\beta\in
\Gamma(T_f^{\perp}M^n)$ satisfy (\ref{eq:gnorm}) then $\Fes_{\va,
\beta}$ given by (\ref{eq:fes}) is  a Combescure transform of $f$
with ${\cal S}_{\va, \beta}$ as the corresponding commuting
Codazzi tensor.

As in the surface case, trivial  Combescure transforms of a given
\ii $f\colon\,M^n\to\R^{N}$ are compositions of $f$ with a
homothety and a translation, whose correspondent commuting Codazzi
tensors are constant multiples of the identity tensor.

We say that an immersion $\Fes\colon\,M^n\to\R^{N}$ is a {\it
Christoffel transform\/} of $f$ if it is a nontrivial {\em
conformal\/} Combescure transform of $f$. In the following result
we classify Euclidean submanifolds of dimension $n\geq 3$ that
admit Christoffel transforms.

\begin{theorem}\po\label{prop:cris}  Let $f\colon\,M^n\to\R^{N}$ be an isometric immersion
that admits a Christoffel transform $\Fes\colon\,M^n\to\R^{N}$. Then $(f,\Fes)$ is a
pair of $2$-isothermic submanifolds. More precisely, there exists a $2$-$CP$-net $(E_+,E_-)$
with respect to which the \sffs of both $f$ and $\Fes$ are adapted. Moreover, if $n\geq 3$
then one of the following holds:
\begin{itemize}
\item[{\em (i)}] for every  $p\in M^n$ there exists a local product representation
$\psi\colon\;M_1\times M_2\to U$ of $(E_+, E_-)$ with $p\in
U\subset M^n$, which is an isometry with respect to a Riemannian
product metric $\<\;,\;\>$ on $M_1\times M_2$, such that $f\circ
\psi=f_1\times f_2$ is an extrinsic  product of \iis with respect
to an orthogonal decomposition $\R^{N}=\R^{N_1}\times \R^{N_2}$.
Moreover,  $\Fes\circ \psi=a((-f_1)\times f_2)\circ \psi+v$ for
some $a\neq 0$ and $v \in \R^{N}$.
\item[{\em (ii)}] for every  $p\in M^n$ there exist a local product representation
$\psi\colon\,I\times N\to U$ of $(E_+, E_-)$ with $p\in U\subset
M^n$, which is an isometry with respect to a warped product metric
$\<\;,\;\>$ on $I\times N$ with warping function $\rho$, an
isometry  $\Phi\colon\;\R_+^{m}\times_{\sigma} \Sf^{N-m}\to \R^N$,
a unit speed curve $\gamma\colon\,I\to \R_+^{m}$ and an \ii
$g\colon\,N\to \Sf^{N-m}$ such that $\gamma_m=\rho=\sigma \circ
\gamma$ and  $f\circ \psi=\Phi\circ(\gamma\times g)$. Moreover,
$\Fes\circ \psi=\Phi\circ(a\tilde{\gamma}\times g)+v$, where
$a\neq 0$, $v\in \R^N$  and $\tilde{\gamma}=\int
\gamma_m^{-2}(\tau)\gamma'(\tau)d\tau$.
\end{itemize}
\end{theorem}
\proof  Since the metric induced by $\Fes$ is
$\<X,Y\>_*=\<d\Fes(X),d\Fes(Y)\>=\<{\cal S} X, {\cal S} Y\>$ for
all $X,Y\in \Gamma(TM^n)$, the symmetry of ${\cal S}$ and the
assumption that $f$ and $\Fes$ are conformal imply that ${\cal
S}^2=\lambda^2I$ for some $\lambda\in C^\infty(M^n)$.  Therefore,
either ${\cal S}=\pm\lambda I$ or $TM^n$ splits as an orthogonal
direct sum $TM^n=E_+\oplus E_-$, where $E_+$ and $E_-$ are the
eigenbundles of ${\cal S}$ correspondent to the eigenvalues
$\lambda$ and $-\lambda$, respectively. In the first case,
$\lambda$ must be a constant $a\neq 0$  by Proposition
\ref{prop:codtensor}-$(i)$, whence $\Fes =af+v$ for some $v\in
\R^{N}$. In the latter case, it follows from Proposition
\ref{thm:codtensor} that $(E_+,E_-)$ is a $CP$-net. Moreover,
since ${\cal S}$ is commuting, the \sff of $f$ is adapted to
$(E_+,E_-)$ and, because of (\ref{eq:rsff}), the same holds for
the \sff of $\Fes$.

  Now assume that $n\geq 3$. Then either $E_+$ or $E_-$, say, the latter, has dimension at
  least two,
  and hence $\lambda$ must be constant along $E_-$ by Proposition \ref{prop:codtensor}-$(i)$.
  Thus Proposition~\ref{prop:codazzi23} applies. In case  $(i)$, since the \sff of $f$
  is adapted to $(E_+,E_-)$, it follows from  the main lemma in \cite{mo} that
  $f\circ \psi=f_1\times f_2$
  splits as an extrinsic  product of \iis with respect to an orthogonal decomposition
  $\R^{N}=\R^{N_1}\times \R^{N_2}$. Moreover, integrating
  $d(\Fes \circ \psi)=d(f \circ \psi)\circ (d\psi)^{-1}\circ{\cal S}\circ d\psi$
  with $(d\psi)^{-1}\circ{\cal S}\circ d\psi=a(\Pi_2-\Pi_1)$ for some $a \neq 0$
  yields $\Fes\circ \psi=a((-f_1)\times f_2)\circ \psi+v$
  for some $v\in \R^{N}$. In case~$(ii)$,  N\"olker's theorem \cite{nol} implies that
there exist an isometry $\Phi\colon\;\R_+^{m}\times_{\sigma}
\Sf^{N-m}\to \R^N$, a unit speed curve $\gamma\colon\,I\to
\R_+^{m}$ and an \ii $g\colon\,N\to \Sf^{N-m}$ such that
$\gamma_m=\rho=\sigma\circ \gamma$ and $f\circ
\psi=\Phi\circ(\gamma\times g)$. Finally, integrating
  $d(\Fes \circ \psi)=d(f \circ \psi)\circ (d\psi)^{-1}\circ{\cal S}\circ d\psi$
  with $(d\psi)^{-1}\circ{\cal S}\circ d\psi=a(\rho\circ \pi_1)^{-2}(\Pi_2-\Pi_1)$
  for some $a\neq 0$ implies
that $\Fes\circ \psi$ is as stated.\qed

 \section[Conformal sphere
congruences]{Conformal sphere congruences}

The  Ribaucour transformation for surfaces in $\R^3$ was extended
as follows to Euclidean submanifolds of arbitrary dimension and
codimension (\cite{dt1}, \cite{dt2}). Two pointwise distinct
immersions $f\colon\, M^n\to\R^{N}$ and
\mbox{$\tilde{f}\colon\,M^n\to \R^{N}$} are said to be related by
a Ribaucour transformation (or each one of them is a  Ribaucour
transform of the other) if  there exist a vector bundle isometry
$\P\colon f^*\Tes\R^{N}\!\to\tilde{f}^*\Tes\R^{N}$, a symmetric
tensor $D$ on $M^n$ with respect to the metric induced by $f$ and
a nowhere vanishing smooth vector field $\delta\in
\Gamma(f^*\Tes\R^{N})$ such that
\begin{itemize}
\item[$(i)$] $\P Z-Z=\<\delta,Z\>(f-\tilde{f})$ for all $ Z\in \Gamma(f^*\Tes\R^{N})$;
\item[$(ii)$] $d\tilde{f}=\P\circ df \circ D$.
\end{itemize}

   Geometrically,  $f$ and $\tilde{f}$ are tangent at each $p\in M^n$ to a
   common $n$-sphere $S(p)$ and $\P$ is the reflection in the hyperplane
   orthogonal to $\tilde{f}(p)-f(p)$. In classical terminology, $f$ and $\tilde{f}$ envelope a
   common $n$-sphere congruence. If $M^n$ is simply-connected, it was shown in
\cite{dt2} (see Theorem $17$) that there exists $(\va,\beta)$
satisfying (\ref{eq:gnorm})  such~that \be\label{eq:rb} \tilde{f}=
f - 2\nu\varphi \Fes, \ee where $\Fes=df(\grad\,\va)+\beta$ and
$\nu^{-1}=\<\Fes,\Fes \>$. Therefore  $\tilde{f}$ is completely
determined by $(\va,\beta)$, or equivalently, by $\va$ and $\Fes$.
We denote \mbox{$\tilde{f}={\cal R}_{\va,\beta}(f)$.} Moreover,
${\cal P}$, $D$ and $\delta$ are given in terms of $(\va,\beta)$
by \be\label{eq:pdo} {\cal P}Z=Z - 2\nu \<\Fes,Z\>\Fes,\;\;\; D=I
- 2\nu\varphi{\cal S}_{\va,\beta}\;\;\;
\mbox{and}\;\;\;\delta=-\varphi^{-1}\Fes. \ee Conversely, given
$(\va,\beta)$ satisfying (\ref{eq:gnorm}) on an open subset
$U\subset M^n$ where  $D$ is invertible, then $\tilde{f}$ given by
(\ref{eq:rb}) defines a \rtf of $f|_U$. The induced metrics,
Levi-Civita connections and \sffs of $f$ and $\tilde{f}$ are
related by \be\label{eq:met} \<X,Y\>^{\sim}=\<DX,DY\>, \ee
\be\label{eq:con3} D\tilde{\nabla}_{X}Y = \nabla_X DY+2\nu \<{\cal
S} X,DY\>\grad\,\va  -2\nu \<\grad\,\va,DY\>{\cal S} X, \ee \be
\label{eq:sff2} \tilde{\a}(X,Y)=\P\left(\a(DX,Y)+2\nu \<{\cal S}
X,DY\>\beta\right). \ee Equation (\ref{eq:sff2}) clarifies the
meaning of the symmetry of the tensor $D$: for each  $\xi\in
T^\perp_fM$ the shape operators $A^f_\xi$ and
$\tilde{A}^{\tilde{f}}_{\P\xi}$ of $f$ and $\tilde{f}$,
respectively,  commute.

\begin{remark}\po\label{re:invar}{\em  For later use we observe the following invariance
property of the Ribaucour transformation. If
\mbox{$\tilde{f}\colon\,M^n\to \R^{N}$}  is a \rtf of
$f\colon\,M^n\to\R^{N}$  with data $(\P,D,\delta)$ and
$\psi\colon\,\tilde{M}^n\to M^n$ is a diffeomorphism,  then the
tensor $\bar{D}$ on $\tilde{M}^n$ defined by $d\psi\circ
\bar{D}=D\circ d\psi$ is symmetric with respect to the metric
induced by $f\circ \psi$, and $\tilde{f}\circ \psi$ is a \rtf of
 $f\circ \psi$ with data $(\bar{\P},\bar{D},\bar{\delta})$, where $\bar{\P}=\P\circ \psi$
and $\bar{\delta}=\delta\circ \psi$. Moreover, if $\tilde{f}={\cal
R}_{\va,\beta}(f)$  then $\tilde{f}\circ \psi={\cal
R}_{\bar{\va},\bar{\beta}}(f\circ \psi)$ for $\bar{\va}=\va\circ
\psi$ and $\bar{\beta}=\beta\circ \psi$. }
\end{remark}

   It follows from (\ref{eq:met}) and the symmetry of $D$ that if $f$ and
   $\tilde{f}$ induce conformal metrics on $M^n$ then $D^2=r^2I$ for some
   $r\in C^\infty(M^n)$. Therefore, either  $D=\pm r I$ or
   $TM^n$ splits  orthogonally as $TM^n=E_+\oplus E_-$, where
$E_+$ and $E_-$ are the eigenbundles of $D$ correspondent to the
eigenvalues $r$ and $-r$, respectively. Since $D=I-2\nu\va {\cal
S}$ by (\ref{eq:pdo}), in the first case ${\cal S}$ must be a
constant  multiple of the identity tensor by Proposition
\ref{prop:codtensor} -$(i)$, in which case the proof of Corollary
$32$ of \cite{dt2} implies that there exists an inversion $I$ in
$\R^{N}$ such that $L'(\tilde{f})=I(L(f))$, where $L$ and $L'$ are
compositions of a homothety and a translation. We say that
$\tilde{f}$ is a {\it Darboux transform \/}  of $f$ if the second
possibility holds, in which case $E_+$ and $E_-$ are also the
eigenbundles of ${\cal S}$ correspondent to its distinct
eigenvalues $\lambda=h(1-r)$ and $\mu=h(1+r)$, respectively, where
$h=(2\nu\va)^{-1}$. Thus, $\tilde{f}$ is a
 Darboux transform  of $f$ if and only if the associated Codazzi tensor ${\cal S}$ has exactly
 two distinct eigenvalues $\lambda, \mu$ everywhere satisfying
\be\label{eq:darboux}(\lambda+\mu)\va=\nu^{-1}=\<\Fes,\Fes\>. \ee

  We now classify  Euclidean submanifolds of dimension $n\geq 3$ that admit Darboux transforms.

\begin{theorem}\po\label{prop:comb2}  Let $f\colon\,M^n\to\R^{N}$ be an isometric immersion
  admitting a Darboux transform $\tilde{f}={\cal
R}_{\va,\beta}(f)\colon\,M^n\to\R^{N}$. Then $(f,\tilde{f})$ is a
pair of $2$-isothermic submanifolds. More precisely, there exists
a $2$-$CP$-net $(E_+,E_-)$ with respect to which the \sffs of both
$f$ and $\tilde{f}$ are adapted. Moreover, if $n\geq 3$ then for
every $p\in M^n$ there exist a local product representation
$\psi\colon\;M_1\times M_2\to U$ of $(E_+, E_-)$ with $p\in
U\subset M^n$, a homothety $H$ and an inversion $I$ in $\R^{N}$
 such that
one of the following holds:
\begin{itemize}
\item[{\em (i)}] $\psi$ is a conformal diffeomorphism with respect to a Riemannian
product metric
on $M_1\times M_2$ and $f\circ \psi=H\circ  I\circ g$, where
 $g=g_1\times g_2\colon\,M_1\times M_2\to\R^{N_1}\times \R^{N_2}=\R^{N}$ is an extrinsic
 product of isometric immersions.
Moreover,  there exists $i\in \{1,2\}$ such that either $M_i$ is
one-dimensional or $g_i(M_i)$ is contained in  some sphere
$\Sf^{N_i-1}(P_i;r_i)\subset \R^{N_i}$.
\item[{\em (ii)}] $\psi$ is a conformal diffeomorphism with respect to a warped product metric
on $M_1\times~M_2$  and $f\circ \psi=H\circ  I\circ
\Phi\,\circ(g_1\times g_2)$, where
$\Phi\colon\;\R_+^{m}\times_{\sigma} \Sf^{N-m}\to \R^N$ is an
isometry and $g_1\colon\,M_1\to \R_+^{m}$ and  $g_2\colon\,M_2\to
\Sf^{N-m}$ are isometric immersions.
\end{itemize}
Conversely, any such \ii admits a Darboux transform.
\end{theorem}
\proof Let $TM=E_+\oplus E_-$ be the orthogonal splitting of $TM$
by the eigenbundles of $D$ correspondent to the eigenvalues $r$
and $-r$, respectively, where $r\in C^\infty(M)$. Since $E_+$ and
$E_-$ are also the eigenbundles of ${\cal S}$ correspondent to the
eigenvalues $\lambda=h(1-r)$ and $\mu=h(1+r)$, respectively, where
$h=(2\nu\va)^{-1}$, it follows from Proposition
\ref{prop:codtensor} that $E_+$ and $E_-$ are umbilical
distributions
 with mean curvature normals
$$\eta^+=(\lambda-\mu)^{-1}(\grad \,\lambda)_{E_-}\;\;\;\;\mbox{and}\;\;\;\;\eta^-=(\mu -\lambda)^{-1}(\grad\, \mu)_{E_+}.$$
We claim that $(E_+,E_-)$ is a $CP$-net. By Proposition
\ref{thm:codtensor},  this is the case if and only if
\be\label{eq:claim}2r X_+(h)X_-(h)+h X_+(h)X_-(r)+h X_-(h)X_+(r)-r
h\, \hess h(X_+,X_-)=0.\ee
 Here and in the sequel, $X_+$ and $X_-$ will always denote sections of $E_+$ and $E_-$,
 respectively.  We now compute $\hess h(X_+,X_-)$. Using (\ref{eq:fes}) we have
$$X(\nu^{-1})=X\<\Fes,\Fes\>=2\<d\Fes(X),\Fes\>=
2\<df({\cal S} X),\Fes\>=2\<{\cal S} X,\grad\,\va\>,$$ whence
$\grad\,\nu=-2\nu^2{\cal S}\, \grad\, \va$ and
$$\grad\, h^{-1}=2(\va\,\grad\, \nu + \nu\,\grad\, \va)=2(-2\nu^2\va\,{\cal S}\,\grad\, \va+
\nu\,\grad\, \va)=2\nu D\,\grad\, \va.$$ It follows that
\be\label{eq:gpsi} \grad\, h=-2h^2\nu D\,\grad\, \va, \ee and thus
$$\nabla_X\grad\, h=-4h X(h)\nu D\,\grad\,\va-2h^2
X(\nu)D\,\grad\,\va-2h^2\nu\nabla_XD\,\grad\, \va.$$ Therefore,
\be\label{eq:hesspsi}
\begin{array}{l}{\displaystyle \hess h(X_+,X_-)= \<\nabla_{X_-}\grad\, h,X_+\>
=-4h X_-(h)\nu \<D\,\grad\,\va,X_+\>-}\vspace{1.5ex}\\
\hspace*{19ex} {\displaystyle 2h^2X_-(\nu)\<D\,\grad\, \va,X_+\>-
2h^2\nu\<\nabla_{X_-}D\,\grad \va,X_+\>.}
\end{array}
\ee We compute  each of the three terms in the right-hand-side of
$(\ref{eq:hesspsi})$. By (\ref{eq:gpsi}), we have that
$X_-(h)=2h^2\nu r X_-(\va)$ and $X_+(h)=-2h^2\nu r X_+(\va)$.
Setting $r=-\va/\tau$, these equations can be rewritten as
$X_-(h)=-(h/\tau)X_-(\va)$ and $X_+(h)=(h/\tau)X_+(\va).$ Using
that $\<D\,\grad\, \va,X_+\>=r X_+(\va)=(r\tau/h)X_+(h)$, we
obtain \be\label{eq:h1}-4h X_-(h)\nu
\<D\,\grad\,\va,X_+\>=-4r\tau\nu
X_-(h)X_+(h)=2h^{-1}X_-(h)X_+(h).\ee Moreover,  $X_-(\nu)=
-2\nu^2\<{\cal S}\,\grad\,
\va,X_-\>=-2\nu^2h(1+r)X_-(\va)=2\nu^2(1+r)\tau X_-(h)$ gives
\be\label{eq:h2}-2h^2X_-(\nu)\<D\,\grad
\va,X_+\>=-4h\nu^2r(1+r)\tau^2
X_-(h)X_+(h)=-\frac{1+r}{rh}X_-(h)X_+(h).\ee It remains to compute
$\<\nabla_{X_-}D\,\grad\, \va,X_+\>.$ By (\ref{eq:con3}) we have
$$D\,\tilde{\nabla}_X Y=\nabla_X DY +2\nu\<{\cal S} X,DY\>\grad\, \va-2\nu\<\grad\, \va,DY\>{\cal S} X.$$
On the other hand, since $\<\;,\;\>^\sim=r^2\<\;,\;\>$, the
connections $\nabla$ and $\tilde{\nabla}$ are also related by
$$\tilde{\nabla}_X Y=\nabla_X Y +r^{-1}(\<\grad\, r,X\>Y+\<\grad\, r,Y\>X-\<X,Y\>\grad\, r).$$
Thus,
$$\begin{array}{l}{\displaystyle  \nabla_X D\,\grad\, \va=D\, \nabla_X \grad\, \va +
r^{-1}\<\grad\, r,X\>D\,\grad\, \va+r^{-1}\<\grad\, r,\grad\, \va\>DX-}\vspace{1.5ex}\\
\hspace*{5ex} {\displaystyle r^{-1}\<X,\grad\, \va\>D\,\grad\,
r-2\nu\<{\cal S} X,D\,\grad\, \va\>\grad\, \va+2\nu\<\grad\,
\va,D\,\grad\, \va\>{\cal S} X.}
\end{array}
$$
We obtain
\begin{eqnarray*}\<\nabla_{X_-}D\,\grad\, \va,X_+\>&=&X _-(r)X_+(\va)-X_-(\va)X_+(r)+2\nu h(1+r)r X_-(\va)X_+(\va)\\
&=&\frac{\tau}{h}(X_-(r)X_+(h)+X_+(r)X_-(h))+\tau(1+r)X_-(h)X_+(h).
\end{eqnarray*}
Hence, \be\label{eq:h3}-2h^2\nu\<\nabla_{X_-}D\,\grad\,
\va,X_+\>=r^{-1}(X_-(r)X_+(h)+X_+(r)X_-(h))+\frac{1+r}{rh}X_-(h)X_+(h).
\ee Then (\ref{eq:claim}) follows by computing  from
(\ref{eq:hesspsi}), (\ref{eq:h1}), (\ref{eq:h2}) and (\ref{eq:h3})
that
$$ \hess h(X_-,X_+)=2h^{-1}X_-(h)X_+(h)+r^{-1}(X_-(r)X_+(h)+X_+
(r)X_-(h)).$$

    Now, since ${\cal S}$ is commuting, the \sff of $f$ is adapted to $(E_+,E_-)$.
    Because of (\ref{eq:sff2}), the same holds for  the \sff of $\tilde{f}$.

     From now on we assume that $n\geq 3$. By Theorem \ref{cor:drh},  for every
     $p\in M^n$ there exists a
     local product representation $\psi\colon\;M_1\times M_2\to U$ of $(E_+, E_-)$
     with $p\in U\subset M^n$,
     which is a conformal diffeomorphism with respect to a Riemannian product metric
     $\<\;,\;\>=\pi_1^*\<\;,\;\>_1+\pi_2^*\<\;,\;\>_2$ on $M_1\times M_2$.
     Therefore, Theorem \ref{thm:moore2} can be applied to $f\circ \psi$. We obtain
     (see Remark \ref{re:main}-$(a)$) that there exist a homothety $H$
      and an inversion $I$ in $\R^{N}$ with respect to a sphere $\Sf^{N-1}(P_0)$
      of unit radius
      such that one of the following holds :
\begin{itemize}
\item[{\em $(a)$}]  $f\circ
\psi=H\circ I\circ (g_1\times g_2)$, where $g_i\colon\;M_i\to
\R^{N_i}$, $1\leq i\leq 2$, are isometric immersions and
$\R^{N}=\R^{N_1}\times \R^{N_2}$ is an orthogonal decomposition.
\item[{\em $(b)$}] after relabelling if necessary,
$f\circ \psi =H\circ I\circ \Phi\circ(f_1\times f_2)$, where
$\Phi\colon\,\R_+^m\times_{\sigma}\Sf^{N-m}\to
 \R^N$ is an isometry,  $f_1\colon\;M_1\to \R_+^m$ is a conformal immersion
with conformal factor $\rho:=(\sigma \circ f_1)$ and
$f_2\colon\;M_2\to \Sf^{N-m}$ is an isometric immersion.
\end{itemize}

       In the latter case,  define $g_1:=f_1\colon\;(M_1,\<\;,\;\>_1^\sim)\to \R_+^m$,
       where $\<\;,\;\>_1^\sim=\rho^2\<\;,\;\>_1$, and set $g_2=f_2$. Then $g_1$ and $g_2$ are
       isometric immersions and $\psi\colon\;(M_1\times M_2,\<\;,\;\>^\sim)\to U$ is still
       a conformal
        diffeomorphism with respect to the warped product metric
        $\<\;,\;\>^\sim=\pi_1^*\<\;,\;\>_1^\sim +(\rho\circ\pi_1)^2\pi_2^*\<\;,\;\>_2.$
        This gives case $(ii)$ of the statement.

        In order to complete the proof of the direct statement, it remains to show that
        the restriction in case $(i)$ holds.
     First notice that, by Remark \ref{re:invar} and the fact that $\psi$ is conformal,
     we have that
        $\tilde{f}\circ \psi$ is a Darboux transform $\tilde{f}\circ \psi={\cal R}_{\bar{\va},\bar{\beta}}(f\circ \psi)$
        of $f\circ \psi$, where $\bar{\va}=\va\circ \psi$ and $\bar{\beta}=\beta\circ \psi$.
        Thus ${\cal S}_{\bar{\va},\bar{\beta}}$ is given by
        $d\psi\circ {\cal S}_{\bar{\va},\bar{\beta}}=
        {\cal S}_{\va,\beta}\circ d\psi$,
        and therefore the eigenbundle net of ${\cal S}_{\bar{\va},\bar{\beta}}$
        is the product net of $M_1\times M_2$. It now follows from
        Proposition~$31$ and equation $(45)$ in \cite{dt2}
        that $\tilde{f}\circ \psi=H\circ I\circ \tilde{g}$,
        where $\tilde{g}={\cal R}_{\hat{\va},\hat{\beta}}(g)$
        is a Darboux transform of $g=g_1\times g_2$
        such that
        ${\cal S}_{\hat{\va},\hat{\beta}}=\tau I+\|g-P_0\|^2{\cal S}_{\bar{\va},\bar{\beta}}$,
         with $\tau=2(\bar{\va}-\<g-P_0,\Fes_{\bar{\va},\bar{\beta}}\>)$. In particular,
        the eigenbundle net of ${\cal S}_{\hat{\va},\hat{\beta}}$ is also the product net of
        $M_1\times M_2$. It now follows easily from Proposition \ref{prop:codtensor} that
either one of the factors is one-dimensional or there exist $a_1,
a_2\in \R$ with $a_1\neq a_2$  such that
        ${\cal S}_{\hat{\va},\hat{\beta}}=a_1\Pi_1+a_2\Pi_2$.     Thus, it suffices to prove
        that in the latter case there must exist $i\in \{1,2\}$
     such that $g_i(M_i)$ is contained in some sphere
     $\Sf^{N_i-1}(P_i;r_i)\subset \R^{N_i}$. Integrating $d{\Fes_{\hat{\va},\hat{\beta}}}=
     dg\circ{\cal S}_{\hat{\va},\hat{\beta}}$ and $d\hat{\va}=
     \<\Fes_{\hat{\va},\hat{\beta}},dg\>$ gives
$$\Fes_{\hat{\va},\hat{\beta}}=a_1(g_1-P_1)\times a_2(g_2-P_2)\;\;\;
\mbox{and}\;\;\;2\hat{\va}=a_1\|g_1\circ
\pi_1-P_1\|^2+a_2\|g_2\circ \pi_2-P_2\|^2+C,$$ for some $P_1\in
\R^{N_1}$, $P_2\in \R^{N_2}$ and $C\in \R$. Using that
$\hat{\va}(a_1+a_2)={\hat{\nu}}^{-1}=
\<\Fes_{\hat{\va},\hat{\beta}},\Fes_{\hat{\va},\hat{\beta}}\>$,
as follows from (\ref{eq:darboux}), we obtain
$$a_1(a_2-a_1)\|g_1\circ \pi_1-P_1\|^2+a_2(a_2-a_1)\|g_2\circ \pi_2-P_2\|^2+C(a_1+a_2)=0,$$
which implies that $a_1\|g_1-P_1\|^2$ and  $a_2\|g_2-P_2\|^2$ must
be constants $a_1r_1^2$
 and $a_2r_2^2$, respectively. Since $a_1$ and $a_2$ can not be both zero,
 we conclude that $\|g_i-P_i\|^2=r_i^2$ for at least one $i\in \{1,2\}$.

    We now prove the converse. It suffices to show that if $g$ is either an extrinsic product
    of
    \iis $g=g_1\times g_2\colon\,M^n=M_1\times M_2\to\R^{N_1}\times
    \R^{N_2}=\R^{N}$ as in $(i)$,
     or it is given by
     $g=\Phi\,\circ(g_1\times g_2)\colon\,M^n=M_1\times_\rho
     M_2\to\R^{N}$, where $\Phi\colon\;\R_+^{m}\times_{\sigma} \Sf^{N-m}\to
     \R^N$ is an isometry and $g_1\colon\;M_1\to \R_+^m$ and
$g_2\colon\;M_2\to \Sf^{N-m}$ are isometric immersions,
 then $g$ admits a Darboux transform $\tilde{g}$.
 For if this is the case and $f=H\circ I\circ g$ for such a $g$,
 then $H\circ I\circ \tilde{g}$ is a Darboux transform of $f$
 by Proposition $31$ in \cite{dt2}.

     Assume first that $g=g_1\times g_2\colon\,M^n=M_1\times M_2\to\R^{N_1}\times \R^{N_2}=\R^{N}$
     is an extrinsic product of \iis with $\|g_2-P_2\|^2=r_2^2$ for some $r_2> 0$,
     $P_2\in \R^{N_2}$.
     Define $\Fes=g_2\circ \pi_2-P_2$ and $\va=r_2^2$. Then $d\va=0=\<\Fes,dg\>$
      and $d\Fes=dg\circ \Pi_2$.
     Therefore  $\va$ and $\Fes$  determine a Ribaucour transform $\tilde{g}$ of $g$ whose
     associated commuting Codazzi tensor is ${\cal S}=\Pi_2$.
     Moreover, (\ref{eq:darboux}) is satisfied, since the eigenvalues of ${\cal S}$ are
     $\lambda=0$ and $\mu=1$. Hence $\tilde{g}$ is in fact
     a Darboux transform of $g$.

   Now suppose that, say, $M_1$ is an open interval $I\subset \R$, so that
   $\alpha:=g_1\colon\;I\to \R^{N_1}$
   is a unit speed curve.   We may assume that the Frenet curvatures $k_1,\ldots,k_{N_1-1}$  of  $\alpha$ are
   nowhere vanishing and consider the linear first order system of ODE's
\be\label{eq:s1}\left\{\begin{array}{l} (i)\;\lambda'=\beta,\;\;\;\;\;(ii)\;\;
\beta'=\lambda+k_1V_2,\vspace{1.5ex}\\
(iii)\;\;V_2'=-k_1\beta+k_2V_3,\;\;\;\;(iv)\;\;V_{N_1}'=-k_{N_1-1}V_{N_1-1},\vspace{1.5ex}\\
(v)\;\;V_j'=-k_{j-1}V_{j-1}+k_jV_{j+1}, \;\;3\leq j\leq
N_1-1,\vspace{1.5ex}
\end{array}\right.
\ee which has the first integral
\be\label{eq:fint}\lambda^2-\lambda'^2-\sum_{j=2}^{N_1}V_j^2=K\in\R.
\ee Let $(\lambda, \beta=\lambda',V_2,\ldots ,V_{N_1})$ be a
solution of (\ref{eq:s1}) with initial conditions chosen so that
the constant $K$ in the right-hand-side of (\ref{eq:fint})
vanishes. Let $\gamma\colon\;I\to \R^{N_1}$ be defined by
$\gamma=\lambda'\alpha'+\sum_{j=2}^{N_1}V_je_j,$
 where $e_1=\alpha',e_2,\ldots, e_{N_1}$ is the Frenet frame of $\alpha$. Using (\ref{eq:s1}) we obtain that
 $\gamma'=\lambda \alpha'$.  Now set
$\Fes=\gamma\circ \pi_1$ and $\va=\lambda\circ \pi_1.$
 Since $\lambda'=\<\gamma,\alpha'\>$, it follows that $d\va=\<\Fes,dg\>$. Moreover,
 $d\Fes=dg\circ ((\lambda\circ \pi_1)\Pi_1)$ whence   $\va$ and $\Fes$  determine a Ribaucour transform
  $\tilde{g}$ of $g$ whose
     associated commuting Codazzi tensor is ${\cal S}=(\lambda\circ \pi_1)\Pi_1$. Furthermore,
 since (\ref{eq:fint})
 holds with $K=0$ and ${\cal S}$ has eigenvalues $\lambda\circ \pi_1$ and $0$, it follows that (\ref{eq:darboux}) is satisfied, for
$$\va(\lambda\circ \pi_1)=(\lambda\circ \pi_1)^2=\|\gamma\circ \pi_1\|^2=\<\Fes,\Fes\>.$$
We conclude that $\tilde{g}$ is  a Darboux transform of $g$.

Finally, we prove that $g=\Phi\,\circ(g_1\times
g_2)\colon\,M^n=M_1\times_\rho M_2\to\R^{N}$,  where
$\Phi\colon\;\R_+^{m}\times_{\sigma} \Sf^{N-m}\to
     \R^N$ is an isometry and $g_1\colon\;M_1\to \R_+^m$ and
$g_2\colon\;M_2\to \Sf^{N-m}\subset \R^{N-m+1}$ are isometric
immersions, also admits a Darboux transform.  Let $g$ be
parameterized by $g=(h_1,\ldots,h_{m-1},h_{m}g_2),$ where
$g_1=(h_1,\ldots,h_{m-1},h_{m})$, $\rho=h_{m}$ and
$g_2=(k_1,\dots,k_{N-m+1})$ has unit length. Define
$$\Fes=(0,\ldots,0,k_1\circ \pi_2,\dots,k_{N-m+1}\circ \pi_2)\;\;\;\mbox{and}\;\;
\va=h_m\circ \pi_1.$$ Then $d\va=d(h_m\circ \pi_1)=\<\Fes, dg\>$
 and $d\Fes=dg\circ {\cal S}$, where
${\cal S}=(h_m\circ \pi_1)^{-1}\Pi_2$. It follows that $\va$ and
$\Fes$  determine a Ribaucour transform $\tilde{g}$ of $g$ whose
associated commuting Codazzi tensor is ${\cal S}$. Moreover, since
the eigenvalues of ${\cal S}$ are $\lambda=0$ and $\mu=(h_m\circ
\pi_1)^{-1}$, we have that $\va(\lambda+\mu)=1=\<\Fes,\Fes\>$.
Thus (\ref{eq:darboux}) is satisfied, whence $\tilde{g}$ is a
Darboux transform of $g$.\qed


 \noindent {\em Acknowledgment.} The
author is grateful to the referee for suggesting a much shorter
proof of Theorem \ref{thm:moore2} than that in the original
manuscript.

{\renewcommand{\baselinestretch}{1}
\hspace*{-25ex}\begin{tabbing}
\indent  \=  Universidade Federal de S\~ao Carlos \\
\indent  \= Via Washington Luiz km 235 \\
\> 13565-905 -- S\~ao Carlos -- Brazil \\
\> e-mail: tojeiro@dm.ufscar.br \\
\end{tabbing}}

\begin{thebibliography}{lbllll}

\bibitem[{\bf Bi}]{bi}  BIANCHI, L.:  {\em Lezioni di Geometria
Differenziale\/}, Bologna, 1927.


\bibitem[{\bf Bo}]{bo} BONNET, O.:  M\'emoire sur la th\'eorie  des
surfaces applicables. {\it J. de l'\'Ecole Polytechnique XLII (1867).}




\bibitem[{\bf Bu}]{bu} BURSTALL, F.: Isothermic surfaces: conformal geometry,
Clifford algebras and integrable systems. Preprint
math-dg/0003096.





\bibitem[{\bf Ce}]{ce1} CECIL, T. E.: {\em Lie Sphere Geometry\/}, Springer--Verlag (1992).



\bibitem[{\bf Ch}]{chris} CHRISTOFFEL, E.:
Ueber einige allgemeine Eigenshaften der Minimumsfl\"achen, {\it
Crelle's J.} {\bf 67} (1867), 218-228.


\bibitem[{\bf DFT}]{dft} DAJCZER, M.; FLORIT, L.; TOJEIRO, R.: On a class of submanifolds
carrying an extrinsic
umbilic foliation, {\it Israel J. of Math.} {\bf 125} (2001),
203-220.



\bibitem[{\bf DT$_1$}]{dt1} DAJCZER, M.; TOJEIRO, R.: An extension of the classical
Ribaucour transformation,
{\it Proc. London Math. Soc.} {\bf 85} (2002), 211--232.

\bibitem[{\bf DT$_2$}]{dt2} DAJCZER, M.;  TOJEIRO, R.: Commuting Codazzi tensors and the
Ribaucour
transformation for submanifolds. {\it Result. Math.} {\bf 44}
(2003), 258--278.




\bibitem[{\bf Da$_1$}]{da1} DARBOUX, G.: Sur les surfaces isothermiques,
{\em C.R.Acad. Sci. Paris} {\bf 128} (1899), 1483--1487.

\bibitem[{\bf Da$_2$}]{da2} DARBOUX, G.: {\it Le\c cons sur la th\'eorie des surfaces\/}
(Reprinted by Chelsea Pub. Co., 1972), Paris 1914.

\bibitem[{\bf H-J}]{hj} HERTRICH-JEROMIN, U.: {\em Introduction to M\"obius differential
geometry\/}, London Math. Lect. Notes Series, vol. $300$,
Cambridge Univ. Press, Cambridge, $2003$.



\bibitem[{\bf La}]{la} LAFONTAINE, J.: Conformal geometry from the Riemannian point of view. Aspects of Mathematics, E 12, Vieweg, Braunschweig, 1988.



\bibitem[{\bf MRS}]{mrs} MEUMERTZHEIM, M.; RECKZIEGEL, H.; SCHAAF. M.: Decomposition of
twisted and warped product
nets,
{\it Result. Math.\/} {\bf 36} (1999), 297--312.


\bibitem[{\bf Mo}]{mo} MOORE, J. D.:  Isometric immersions of Riemannian  products,
{\it J. Diff. Geom.\/}~{\bf 5} (1971), 159--168.

\bibitem[{\bf No}]{nol} N\"OLKER, S.:  Isometric immersions of warped products,
{\it Diff. Geom. Appl.\/}~{\bf 6} (1996), 31--50.


\bibitem[{\bf RS}]{rs} RECKZIEGEL, H.; SCHAAF. M.: De Rham decomposition of netted
manifolds,
{\it Result. Math.\/} {\bf 35} (1999), 175--191.

\bibitem[{\bf Re}]{re} RECKZIEGEL, H.: Kr\"ummungsflachen von isometrischen Immersionen in
R\"aumen konstanter Kr\"ummung, {\it Math. Ann.\/} {\bf 223}
(1976), 169--181.



\bibitem[{\bf To}]{to} TOJEIRO, R.:   Conformal De Rham decomposition of Riemannian
manifolds. To appear in {\it Houston J. Math}.


\end{thebibliography}
\end{document}